# An Ultra-high-speed Reproducing Kernel Particle Method

Siavash Jafarzadeh[1,3,✉], Michael Hillman[2,3]

[1] *Department of Mathematics, Lehigh University, Bethlehem, PA 18015, USA*

[2] *Karagozian & Case, 700 N brand Blvd, # 700, Glendale, CA 91203, USA*

[3] *Department of Civil and Environmental Engineering, The Pennsylvania State University, University Park, PA 16802, USA*

**Abstract**

In this work, the fast-convolving reproducing kernel particle method (FC-RKPM) is introduced. This method is hundreds to millions of times faster than the traditional RKPM for 3D meshfree simulations. In this approach, the meshfree discretizations with RK approximation are expressed in terms of convolution sums. Fast Fourier transform (FFT) is then used to efficiently compute the convolutions. Certain modifications to the domain and shape functions are considered to maintain generality for complex geometries and arbitrary boundary conditions. The new method does not need to identify, store, and loop over the neighbors which is one of the bottleneck of the traditional meshfree methods. As a result, the run-times and memory allocations are independent of the number of neighbors and the shape function's support size. As a model problem, the method is laid out for a Galerkin weak form of the Poisson problem with the RK approximation, and is verified in 1D, 2D, and 3D. Tables with run-times and allocated memory are presented to compare the performance of FC-RKPM with the traditional method in 3D. The performance is studied for various node numbers, support size, and approximation degree. All the implementation details and the roadmap for software development are also provided. Application of the new method to nonlinear and explicit problems are briefly discussed as well.

1. ## Introduction

Meshfree methods have been developed over the past decades to overcome the difficulties one faces when using the traditional finite element method, especially when solving problems that involve major topological changes such as large deformation, flow, fracture, etc. [1–6]. Development of kernel-based meshfree methods can be traced back to smoothed particle hydrodynamics (SPH) [7, 8]. Reproducing kernel particle methods (RKPM) [2, 9] are an advanced type of meshfree methods, where a correction function is added to the kernel in the SPH formulation, such that approximated functions meet certain polynomial reproducing conditions. This leads to significant improvement of the accuracy and convergence of the approximation and corrects for errors near the boundaries. The computational cost of the meshfree methods is however expensive, especially in 3D. Most of the efficient implementations of RKPM require neighbor search and storage, and frequent looping over neighbors for each node. When a weak-form-based RKPM is considered, quadrature is performed that entails a loop over neighbors for strain calculations, and another for assembly. The RK shape functions themselves are also expensive to construct since they require several vector/matrix operations of size four for linear basis (in 3D)[2, 10, 11].

There have been a number of attempts to reduce the computational cost of RKPM. Implicit gradients were introduced as cheap alternatives to the exact derivatives of RK shape functions [12, 13]. Consistent

---

pseudo-derivatives were introduced in [14] in the same spirit. In [15], the authors demonstrated that the correspondence-based PD deformation gradients can also be expressed as the interpolation of relative nodal displacement with an approximation related to the implicit gradient RK shape function, but with a 3x3 instead of a 4x4 correction matrix. In [16–18], the concept of arbitrary order recursive gradients of meshfree shape function was proposed, where the higher order gradients of shape function is simplified by the recursive operation of only the first order gradient of shape function, which significantly reduced the computational cost. To improve efficiency of numerical integration in weak forms, Direct nodal integration (DNI) is the popular choice where nodes are used as the integration points for quadrature [19]. The tradeoff for the efficiency gains from DNI are zero-energy mode instabilities and non-convergence of the solution. The issue then is the efficiency of techniques used to stabilize and correct the solution. In [11] an accelerated stabilization using implicit gradients was introduced to this end. Other attempts to reduce the cost of RKPM include efficient neighbor-search techniques such as partitioning, tree search algorithms, etc. [20, 21], and also vectorizing the moment matrix inversion [22]. While all these techniques result in some computational efficiency, the true bottleneck of meshfree methods remains untouched, which is the need for looping over neighbors for every node. If $N$ is the total number of nodes and $M$ is the number of neighbors within the support of the kernel function, the traditional RKPM has the complexity of $O(NM)$ to compute most terms of the system. Computing and construction of the stiffness matrix, if needed, is even worse and has a complexity of $O(NM^2)$ since a second loop over the neighbors is required.

Interestingly, the computational aspects of the meshfree methods are very similar to those of the nonlocal models, particularly Peridynamics (PD) [23–25]. Peridynamics [26–29] is a nonlocal extension of continuum mechanics where spatial derivatives are replaced with volume integrals, allowing for emergence and evolution of discontinuities in the field variables, making the theory suited for fracture, damage, and fragmentation (e.g., see [3, 4, 30–32]). Similarities between RKPM and the meshfree discretization of a class of peridynamic formulation have been noticed and investigated in several studies [15, 23–25]. As in the RKPM discretization, PD simulations are also computationally expensive, especially in 2D and 3D due to the interactions of nodes within a finite-size proximity and the quadrature over such neighborhood for each node [33]. Recently a fast convolution-based method (FCBM) for peridynamics has been introduced where the convolutional structures of PD volume integrals are exploited to allow for computing quadratures via fast Fourier transform [34–40]. In FCBM, explicit identification of neighbors and looping over them for integration are avoided by transforming the nature of operations from direct summation over the neighbors in the physical space (the common approach), to cheap multiplications in the Fourier space. This is achieved by using convolution theorem and the *fast Fourier transform* (FFT) algorithm at the cost of $O(N\log_2 N)$. While most Fourier and FFT-based methods require periodicity of the domain which limit their applicability to very especial problems, FCBM uses an embedded constraint (EC) approach to enforce the desired boundary conditions on arbitrary geometries, and therefore remains applicable to real world engineering problems.

Inspired by the FCBM developed for PD, in this study, we explore the potentials that FFT-accelerated summation offers within the RK methods. In particular, we introduce the Fast-Convolving RKPM (FC-RKPM), where the summations of the RK approximation of functions and those arising from direct nodal integration in the Galerkin weak form are efficiently computed using FFT operations. The introduced method does not require neighbor search and storage a priori, nor does it require construction and storage of stiffness and/or mass matrices. We verify the method against analytical solutions and provide run-times and memory allocation comparison between the new and the traditional implementations. This work also aims to provide a roadmap for developing softwares based on FC-RKPM. For that reason, all implementation steps are provided along with detailed pseudocodes for various scenarios are presented.

This paper is organized as follows: In Section 2, the fundamentals of RKPM and its traditional implementation are discussed. In Section 3, necessary background information on FFT and fast convolution is provided. The new method is presented in Section 4, and is followed by the implementation details and algorithms in Section 5. Numerical examples for verification and performance comparisons are provided in Section 6.

## 2. Reproducing Kernel Particle Method

In this part, we briefly discuss the reproducing kernel (RK) approximation and the Galerkin weak form used in RKPM. We also describe the method's standard implementation. In this paper, boldface denotes vectors, tensors, and arrays, whereas scalars are denoted by plain letters. For clarity and conciseness, formulations are presented for the 2D case throughout this article. The 1D and the 3D formulations can be easily obtained by dropping or adding one dimension to the equations.

### 2.1. Reproducing Kernel approximation

The core building block of all RK methods is the approximation of functions via RK shape functions. Let the spatial domain of interest $\Omega$ be discretized into a finite set of nodes $S = \{\boldsymbol{x}_1, \boldsymbol{x}_2, \ldots, \boldsymbol{x}_N\}$, where $\boldsymbol{x}_I$ ($I = 1,2,\ldots,N$) denote spatial coordinate of the node $I$; for 2D Cartesian coordinate system, $\boldsymbol{x}_I = \{x_I, y_I\}$. The RK approximation of a function $u(\boldsymbol{x})$ is given by[9, 41]:

$$u^h(\boldsymbol{x}) = \sum_{I=1}^{N} \Psi_I(\boldsymbol{x}) d_I \tag{1}$$

where $\Psi_I(\boldsymbol{x})$ is the *RK shape function* of node $\boldsymbol{x}_I$ evaluated at $\boldsymbol{x}$, and $d_I$ is the corresponding coefficient. In practice the above loop is performed only for non-zero neighbors. The discrete RK shape functions are defined as:

$$\Psi_I(\boldsymbol{x}) = \boldsymbol{H}^{\mathrm{T}}(\boldsymbol{0}) \boldsymbol{M}^{-1}(\boldsymbol{x}) \boldsymbol{H}(\boldsymbol{x} - \boldsymbol{x}_I) \phi_a(\boldsymbol{x} - \boldsymbol{x}_I) \tag{2}$$

where $\boldsymbol{M}(\boldsymbol{x})$ is the *moment matrix*:

$$\boldsymbol{M}(\boldsymbol{x}) = \sum_{I=1}^{N} \boldsymbol{H}(\boldsymbol{x} - \boldsymbol{x}_I) \boldsymbol{H}^{\mathrm{T}}(\boldsymbol{x} - \boldsymbol{x}_I) \phi_a(\boldsymbol{x} - \boldsymbol{x}_I), \tag{3}$$

$\boldsymbol{H}$ is the *vector of monomial basis functions*:

$$\boldsymbol{H}(\boldsymbol{x}) = [1, x, y, x^2, xy, y^2, \ldots, y^n]^{\mathrm{T}}, \text{ and } \boldsymbol{x} = \begin{Bmatrix} x \\ y \end{Bmatrix}. \tag{4}$$

$n$ is the monomial degree to which the approximation is able to reproduce. Given the spatial dimension of the problem ($d$), the minimum number of non-colinear neighbors (non-coplanar in 3D) to ensure an invertible moment matrix $\boldsymbol{M}$ is given by [42]:

$$\frac{(n+d)!}{d! \, n!} \tag{5}$$

$\phi_a$ in Eqs. (2) and (3), is the *Kernel function* with a finite size support. One example of the kernel, which we use in this study, is the function with rectangular support and cubic B-Spline profile in each dimension [1]:

$$\phi_a(\boldsymbol{x}) = \phi_{a_x}(x) \phi_{a_y}(y) \tag{6}$$

$$\phi_{a_x}(x) = \begin{cases} \frac{2}{3} - 4z^2 + 4z^3 & \text{for } 0 \le z \le \frac{1}{2} \\ \frac{4}{3} - 4z + 4z^2 - \frac{4}{3}z^3 & \text{for } \frac{1}{2} \le z \le 1, \quad z = \frac{|x|}{a_x} \\ 0 & \text{for } z > 1 \end{cases} \quad (7)$$

and $\phi_{a_y}(y)$ is defined similarly, where $x$ is replaced with $y$ in Eq. (7). $a_x$ and $a_y$ are respectively the support sizes in $x$ and $y$ coordinate directions. The kernel function dictates the locality of the RK approximation, which can be controlled by the support parameters $a_x$ and $a_y$ as seen from the above.

### 2.2. The Galerkin weak form

The RKPM uses the Galerkin framework and RK approximation to solve partial differential equations (PDEs). For demonstration, consider the following Poisson problem with Dirichlet and Neumann boundary conditions:

$$\begin{cases} \nabla^2 u + r = 0, & \text{on } \Omega \\ u = g, & \text{on } \Gamma_g, \\ \nabla u \cdot \boldsymbol{n} = q, & \text{on } \Gamma_q \end{cases} \quad (8)$$

where $u(\boldsymbol{x})$ is the field variable to find, $r(\boldsymbol{x})$ is a source term, $\Omega$ is the domain of interest, $\Gamma_g$ and $\Gamma_q$ are the subsets of $\Omega$'s boundary, respectively corresponding to the Dirichlet and Neumann boundary conditions (see Figure 1). $\nabla^2 = \boldsymbol{\nabla} \cdot \boldsymbol{\nabla}$ is the Laplacian and $\boldsymbol{\nabla}$ denotes the gradient operator. $\boldsymbol{n}$ is the outward normal unit vector on the boundary, and $u_g(\boldsymbol{x})$ and $q(\boldsymbol{x})$ are known functions.

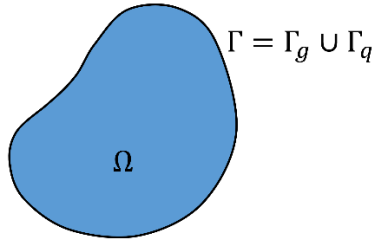

Figure 1. A generic bounded domain ($\Omega$), and its boundary ($\Gamma = \Gamma_g \cup \Gamma_q$).

The Galerkin weak form of the Poisson problem yields:

$$\int_\Omega \boldsymbol{\nabla} w^h \cdot \boldsymbol{\nabla} u^h \, d\Omega - \int_\Omega w^h r \, d\Omega - \int_{\Gamma_q} w^h q \, d\Gamma = 0 \quad (9)$$

Where $u^h$ is the approximated field variable given by Eq. (1), and $w^h(\boldsymbol{x})$ is the *test function* given by:

$$w^h(\boldsymbol{x}) = \sum_{I=1}^N \Psi_I(\boldsymbol{x}) c_I \quad (10)$$

Therefore:

$$\boldsymbol{\nabla} u^h(\boldsymbol{x}) = \sum_{I=1}^N \boldsymbol{\nabla} \Psi_I(\boldsymbol{x}) d_I \quad (11)$$

$$\boldsymbol{\nabla} w^h(\boldsymbol{x}) = \sum_{I=1}^N \boldsymbol{\nabla} \Psi_I(\boldsymbol{x}) c_I \quad (12)$$

where:

$$\nabla\Psi_I(x) = \begin{Bmatrix} \Psi_{I,x} \\ \Psi_{I,y} \end{Bmatrix} \quad (13)$$

Since computing the shape function's exact derivatives are expensive, *implicit gradients* can be used as cheap alternatives which are approximations, obtained by enforcing the reproducing conditions to the derivatives themselves [43]. Approximation with Implicit Gradient allows us to obtain the approximation to partial derivatives of RK shape functions by simply replacing $H^T(0)$ in Eq. (2) with another constant vector that depends on the differentiation order and coordinate direction.

The gradient of RK shape functions in Eqs. (11) and (12) can then be approximated as:

$$\nabla\Psi_I(x) \cong \Psi_I^\nabla(x) = \begin{Bmatrix} \Psi_I^{\nabla x} \\ \Psi_I^{\nabla y} \end{Bmatrix} \quad (14)$$

where $\Psi_I^\nabla(x)$ denotes the Implicit Gradient and:

$$\Psi_I^{\nabla x}(x) = [H^{\nabla x}]^T M^{-1}(x) H(x - x_I) \phi_a(x - x_I) \quad (15)$$

$$\Psi_I^{\nabla y}(x) = [H^{\nabla y}]^T M^{-1}(x) H(x - x_I) \phi_a(x - x_I)$$

with

$$H^{\nabla x} = [0, -1, 0, \dots, 0]^T \quad (16)$$

$$H^{\nabla y} = [0, 0, -1, 0, \dots, 0]^T$$

Note that the size of $H^{\nabla x}$ and $H^{\nabla y}$ is the same as $H$. Details of the Implicit Gradient approximation and finding the appropriate substitutive vectors corresponding to the differentiation is given in [1]. Eqs. (15) and (16) suffice for the examples in this work.

### 2.3. Matrix Forms and Implementation

For computer implementation of RKPM, the following data structure is usually adopted. In this article, nonitalic boldface letters are used to denote the *implementation arrays* which their size depends on the degrees of freedom ($N$). Implementation arrays should not be confused with italic boldface letters used for continuum arrays such as $M$ and $H$ which denote quantities with dimensions independent of discretization size. One can express Eqs. (1), (10), (11) and (12) as:

$$u^h = \mathbf{Nd} \quad (17)$$

$$w^h = \mathbf{Nc} \quad (18)$$

$$\nabla u^h = \mathbf{Bd} \quad (19)$$

$$\nabla w^h = \mathbf{Bc} \quad (20)$$

where

$$\mathbf{d} = [d_1, d_2, \dots, d_N]^T \quad (21)$$

$$\mathbf{c} = [c_1, c_2, \dots, c_N]^T \quad (22)$$

$$\mathbf{N} = [\Psi_1, \Psi_2, \dots, \Psi_N] \quad (23)$$

$$\mathbf{B} = \begin{bmatrix} \Psi_{1,x} & \Psi_{2,x} & \cdots & \Psi_{N,x} \\ \Psi_{1,y} & \Psi_{2,y} & \cdots & \Psi_{N,y} \end{bmatrix} \tag{24}$$

Note that the RK partial derivatives in Eq. (24) can be approximated by Eq.(15).

Substituting the matrix forms above in the Galerkin weak form given by Eq. (9) yields:

$$\mathbf{Kd} = \mathbf{f}^r + \mathbf{f}^q \tag{25}$$

where

$$\mathbf{K} = \int_\Omega \mathbf{B}^T \mathbf{B} \, d\Omega \tag{26}$$

$$\mathbf{f}^r = \int_\Omega \mathbf{N}^T r \, d\Omega \tag{27}$$

$$\mathbf{f}^q = \int_{\Gamma_q} \mathbf{N}^T q \, d\Gamma \tag{28}$$

As briefly stated in the introduction, the computational complexity associated with constructing the stiffness matrix $\mathbf{K}$ is at least of $O(NM^2)$ where $M$ is the number of neighbors, i.e., nodes within the compact support of the kernel. For each node, a loop is needed over the neighbors, and for each pair of nodes, another loop over their common neighbors is required, resulting in a triple nested For-loops with the complexity of $O(NM^2)$. The computational complexity of $\mathbf{f}^r$ and the $\mathbf{Kd}$ product (e.g., needed for Krylov subspace methods) is of $O(NM)$, because they are computed via double nested For-loops.

The cost for boundary terms is different. The computational complexity of $\mathbf{f}^q$ is of $O(N^\Gamma M^\Gamma)$, where $N^\Gamma$ denotes the number of nodes near or on the boundary $\Gamma$ (nodes that their shape function cover at least one node on the boundary), and $M^\Gamma$ is the number of neighbors, for such node, that locate on the boundary. Similarly, if essential boundary conditions are weakly enforced (see [44]), additional terms will appear in Eq. (25), which would have the complexity of $O(N^\Gamma M^\Gamma)$, since they are, too, computed as boundary integrals.

*Remark*: Given that $N^\Gamma$ and $M^\Gamma$ are in a lower dimensional manifold compared to $N$ and $M$, the major computational cost of implicit RKPM analysis is associated with construction of the matrix $\mathbf{K}$ in the first place, and computing $\mathbf{f}^r$ and $\mathbf{Kd}$ product in the second place.

### 2.4. Dirichlet Boundary Conditions

There are several ways to enforce the Dirichlet (essential) boundary conditions in RKPM. The simplest approach is to assume $d_I \cong u(\mathbf{x}_I)$, and strongly set the coefficient values when solving the system given by Eq. (25):

$$d_I = g(\mathbf{x}_I) \text{ for } \mathbf{x}_I \in \Gamma_g, \tag{29}$$

However, since RK shape functions generally lack the Kronecker $\delta$ property, the coefficients are not exactly equal to the values of the field variable, hence Eq. (29) may introduce some error to the solution depending on deviation of $d_I$ from $u(\mathbf{x}_I)$. Several alternative methods and corrective techniques have been proposed for strong enforcement of Dirichlet BCs, e.g., modification of the shape functions near the boundaries [45], transformation method [46], as well as collocating essential boundary conditions in the strong-form-based collocation methods [47]. It is also shown in a study that strong enforcement of essential BC should follow *consistent weak forms*, otherwise the rate of convergence would be negatively affected [44].

Another approach for enforcing Dirichlet BC is to weakly enforce the BC using methods such as penalty, Lagrange multiplier, and Nitsche methods [48]. Note that weak enforcements of essential BC are accompanied by additional terms in the weak forms, and therefore, additional forces in Eq. (27) which would be associated with the BCs. Depending on the weak enforcement method, extra degrees of freedom or penalty parameters may be required.

Since the purpose of this study is deriving the fast-convolving RKPM formulation, its verification, and its performance compared to the standard method, in the presented numerical examples of Section 6, we adopt the simplest and most-practiced approach for essential boundary conditions, which is the strong enforcement shown by Eq. (29). The methods in this paper can be extended to any number of techniques that consistently enforce essential boundary conditions via modifying the approximation, or by employing the techniques for boundary integrals in weak-based methods.

### 3. Preliminaries and background on fast convolution sums

In this part, we briefly describe the basics of discrete Fourier transform, circular convolution, and the fast convolution summation. Again, we provide the equations for 2D; the 1D and the 3D cases can easily follow by dropping or adding one dimension.

*Discrete Fourier Transform*

Let $a(i,j) = \{a_{ij} | i,j = 1,2,...\}$ be a discrete 2D periodic sequence with $a(i,j) = a(i + t_1 N_1, j + t_2 N_2)$, where $N_1$ and $N_2$ are the positive integers that define the period range for the first and the second indices (representing the two dimensions), and $t_1$ and $t_2$ can be any positive integers. The *discrete Fourier transform* (DFT) and its inverse (iDFT) operations can respectively be defined as:

$$\hat{a}_{k_1 k_2} = \sum_{j=1}^{N_2} \sum_{i=1}^{N_1} a_{ij} e^{-2\pi \zeta \left[\frac{(k_1-1)(i-1)}{N_1} + \frac{(k_2-1)(j-1)}{N_2}\right]}, \tag{30}$$

$$a_{ij} = \frac{1}{N_1 N_2} \sum_{k_2=1}^{N_2} \sum_{k_1=1}^{N_1} \hat{a}_{k_1 k_2} e^{2\pi \zeta \left[\frac{(k_1-1)(i-1)}{N_1} + \frac{(k_2-1)(j-1)}{N_2}\right]}, \tag{31}$$

where $\zeta = \sqrt{-1}$, and $\hat{a}_{k_1 k_2} = \hat{a}(k_1, k_2)$ is the discrete Fourier coefficient of the sequence $a$, associated with the mode $K = \{k_1, k_2\}$. Note that the DFT definitions can slightly vary in literature; Eqs. (*30*) and (*31*) are consistent with the version presented in MATLAB's documentation [49].

*Fast Fourier Transform (FFT)*

FFT and its inverse (iFFT), are efficient algorithms to compute DFT and iDFT operations [50, 51]. The computational complexity of FFT and iFFT for multi-dimensional arrays of total $N$ entries is of $O(N\log_2 N)$.

*Discrete circular convolution*

Let $a$ and $b$ be two 2D periodic sequences. The discrete circular convolution can be defined as the periodic sequence $c_{ij}$:

$$c_{ij} = [a \circledast b]_{ij} = \sum_{q=1}^{N_2} \sum_{p=1}^{N_1} a(p,q)\, b(i-p, j-q) \tag{32}$$

*Convolution Theorem:* Let $a$ and $b$ be two periodic sequences (here, of two dimension). Then:

$$\hat{c}_{k_1 k_2} = \hat{a}_{k_1 k_2} \hat{b}_{k_1 k_2} \tag{33}$$

*Fast convolution sums:*

Let **a** and **b** be $d$-dimensional arrays of total $N = N_1 N_2 \ldots N_d$ entries, containing one period of the periodic $d$-dimensional sequences $a$ and $b$. The fast circular convolution sum for **a** and **b** is expressed as:

$$\mathbf{c} = \mathbf{F}^{-1}\{\mathbf{F}(\mathbf{a}) \circ \mathbf{F}(\mathbf{b})\} \tag{34}$$

where the operator ( $\circ$ ) denotes the *elementwise (Hadamard) product* of two arrays, and $\mathbf{F}$ and $\mathbf{F}^{-1}$ denote the *FFT* and its inverse (*iFFT*) operations. The fast convolution sum computes the circular convolution at the cost of $O(N\log_2 N)$, while the direct summation in Eq. (*32*) costs $O(N^2)$.

For linear convolutions (as opposed to circular), where the sequences are not periodic, fast convolution by FFT can be achieved by certain techniques and modifications, such as zero-padding of the convolving sequences to remove the influence of "wrap-arounds" in the circular convolution [52].

*Convolution sums for discrete fields*

Let $\mathbb{T} = [0, L_1] \times [0, L_2]$, be a periodic continuous domain (0 is identified with $L_1$ and $L_2$ in $x$ and $y$ directions due to periodicity) and be uniformly discretized with a finite number of nodes, equally spaced in each coordinate directions:

$$\mathbf{x}_{ij} = \begin{Bmatrix} x_i \\ y_j \end{Bmatrix} = \begin{Bmatrix} (i-1)\Delta x \\ (j-1)\Delta y \end{Bmatrix}; \Delta x = \frac{L_1}{N_1}, \Delta y = \frac{L_2}{N_2}; \text{ and } i = 1,2,\ldots,N_1 \quad j = 1,2,\ldots,N_2 \tag{35}$$

let $a(\mathbf{x})$ and $b(\mathbf{x})$ be two function on $\mathbb{T}$. The discrete circular convolution of $a$ and $b$ is defined as:

$$[a \circledast b]_{ij} = \sum_{q=1}^{N_2} \sum_{p=1}^{N_1} a(x_i, y_j)\, b(x_i - x_p, y_j - y_q) \tag{36}$$

Note that uniform discretization is essential for the discrete circular convolution theorem in Eq. (33) to be directly applicable to Eq. (36). The equivalency between Eq. (*32*) and Eq. (36) is shown below for a uniformly discretized domain:

$$\sum_{q=1}^{N_2}\sum_{p=1}^{N_1} a(x_i, y_j)\, b(x_i - x_p, y_j - y_q) = \sum_{q=1}^{N_2}\sum_{p=1}^{N_1} a(i,j)\, b((i-p)\Delta x, (j-q)\Delta y) \tag{37}$$

$$= \sum_{q=1}^{N_2}\sum_{p=1}^{N_1} a(p,q)\, b(i-p, j-q)$$

*Remark:* For a nonuniform discretization $\Delta x$ and $\Delta y$ are not constant, and Eq. (35) and consequently Eq. (37) do not hold. However, fast convolution sums for non-equispaced discrete field (non-uniform discretization of space), can be carried out by non-uniform FFT [53–56], which are algorithms to compute non-uniform DFT [53] at the cost of $O(N\log_2 N)$. This feature, however, is not explored in this work.

*Notations and arrays for multi-dimensional analysis:*

In this study, two types of indexing are used for representing arrays' elements which contain discrete data for field quantities: one is the *total indexing* where a single identifier is used for numbering the elements, regardless of the spatial dimension (e.g., $\boldsymbol{x}_I$). This notation is concise and best to use when explicit expression of the dimensions is not needed for computation. Like most articles in RKPM, we have used this notation in Section 3 of this study. The other indexing approach is the *dimensional indexing* (e.g., $\boldsymbol{x}_{ij}$ or $\boldsymbol{x}(i,j)$ in 2D) which uses distinct identifiers for each coordinate directions. This indexing is used when explicit expression of the dimensionality is needed in computations, i.e., it is necessary to work with $d$-dimensional arrays to carry out the computations ($d$ being the spatial dimension). As noted in this section, multi-dimensional FFT is a case where computations are carried out using $d$-dimensional arrays, and therefor, dimensional indexing is preferred.

Eq. (38) can be used to translate between the total indexing and dimensional indexing as needed.

$$a_I \ (I = 1, 2, \dots N) \overset{\text{for 2D}}{\Longleftrightarrow} a_{ij}(i = 1, \dots, N_1 \ ; j = 1, \dots, N_2), \quad \text{and } N = N_1 N_2 \tag{38}$$

$$\sum_{I=1}^{N} (\cdot) \overset{\text{for 2D}}{\Longleftrightarrow} \sum_{j=1}^{N_2} \sum_{i=1}^{N_1} (\cdot)$$

The convolution sum of two functions described with multi-dimensional indexing in Eq. (*36*), can be expressed using the total indexing as well:

$$[a \circledast b]_I = \sum_{J=1}^{N} a(\boldsymbol{x}_J) b(\boldsymbol{x}_I - \boldsymbol{x}_J) = \sum_{J=1}^{N} a_J b_{I-J} \tag{39}$$

## 4. The Fast-Convolving Reproducing Kernel Particle Method

In this section, we derive the FC-RKPM formulation and lay out its implementation. The core idea in this method is to express RKPM formulation in terms of convolution sums, and then compute those summations via FFT, following the convolution theorem. It can be easily shown that the standard RKPM formulation described in Section 3, can be expressed in convolutional forms. However, the fast circular convolution would then be only applicable to periodic domains, i.e., problems with periodic BCs. This would limit the method's applicability as most real-life engineering problems are not periodic, rather they are defined on bounded domains with boundary conditions. FC-RKPM for non-periodic cases can be achieved if one uses techniques such as zero-padding with FFT [52] for fast evaluation of the linear convolutions in the original RKPM. In this study, however, inspired by fast convolution-based peridynamics [36, 37], we adopt another approach. We modify the RKPM formulations such that the fast circular convolution becomes applicable to problems over bounded domains with boundary conditions. The original formulation of RKPM then becomes a special case of the modified version. With the proposed approach one does not need to apply techniques (e.g., zero-padding, etc.) for every convolution, as one would need with the original form.

### 4.1. Modifications of the RK shape functions

To achieve the FC-RKPM discretization for problems with general boundary conditions (non-periodic), the following modifications are considered:

First, the bounded domain of interest $\Omega$ for the given problem is extended to a periodic box $\mathbb{T}$ (see Figure 2).

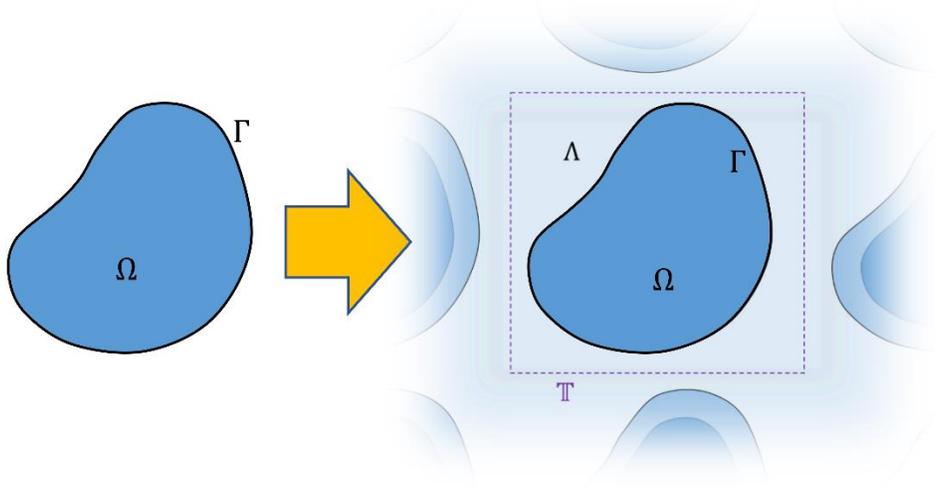

Figure 2. Embedding a generic 2D bounded domain $\Omega$ in a periodic box $\mathbb{T}$.

Then the extended domain $\mathbb{T}$ is uniformly discretized into a total of $N$ nodes, with the position vectors $x_I$ ($I = 1,2,\ldots,N$), that are equally spaced in each coordinate directions.

Next, a characteristic function is used to modify the RK shape functions, such that the solution on the bounded domain is identically represented on the extended domain $\mathbb{T}$. Let

$$\chi(x_I) = \begin{cases} 1 & x_I \in \Omega \\ 0 & x_I \in \Lambda \ (\mathbb{T}\backslash\Omega) \end{cases} \tag{40}$$

be the *domain's characteristic function*. We modify RK shape functions, the Implicit Gradients, and the moment matrices as follows:

$$\Psi_I(x) = \chi(x)\chi(x_I)[b^0]^\mathrm{T}(x)H^a(x - x_I) \tag{41}$$

$$\Psi_I^{\nabla_x}(x) = \chi(x)\chi(x_I)[b^x]^\mathrm{T}(x)H^a(x - x_I) \tag{42}$$

$$\Psi_I^{\nabla_y}(x) = \chi(x)\chi(x_I)[b^y]^\mathrm{T}(x)H^a(x - x_I) \tag{43}$$

where:

$$H^a(x) = H(x)\phi_a(x) \tag{44}$$

$$[b^0]^\mathrm{T}(x) = H^\mathrm{T}(0)M^{-1}(x) = 1^{\text{st}} \text{ row of } M^{-1} \tag{45}$$

$$[b^x]^\mathrm{T}(x) = [H^{\nabla_x}]^T M^{-1}(x) = -2^{\text{nd}} \text{ row of } M^{-1} \tag{46}$$

$$[b^y]^\mathrm{T}(x) = [H^{\nabla_y}]^T M^{-1}(x) = -3^{\text{rd}} \text{ row of } M^{-1} \tag{47}$$

are defined for conciseness of derivations that follows, and

$$M(x) = [1 - \chi(x)]\mathbb{I} + \chi(x) \sum_{I=1}^{N} \chi(x_I) H(x - x_I)[H^a]^{\mathrm{T}}(x - x_I), \quad (48)$$

where $\mathbb{I}$ is the $s \times s$ identity matrix ($s$ is the size of $H$).

Here we explain how using the modified forms on the periodic domain $\mathbb{T}$, are equivalent to using the original version on the bounded domain $\Omega$. Consider the three possibilities below:

- When $x$ and all of its neighbors $x_I$ locates inside $\Omega$, i.e., $\chi(x) = \chi(x_I) = 1$, Eqs. (41) to (43) and Eq. (48) become identical to Eqs. (2), (15) and (3), respectively.
- When $x$ locates inside $\Omega$ but near the boundaries where one or more neighbor, $x_I$, locate in $\Lambda$, i.e., $\chi(x) = 1$ and $\chi(x_I) = 0$, the contributions from nodes outside $\Omega$ to the summation in Eq. (48) will be zero, which is equivalent to the $M$ matrix resulting from Eq. (3) for a bounded domain near the boundaries with a cut-off neighborhood. The shape function of $x_I$ evaluated at $x$ will be also zero: $\Psi_I(x) = 0$ (same for the implicit gradients). This is consistent with what Eqs. (1) and (11) give for approximating $u^h$ and $\nabla u^h$, since such $x_I$ outside $\Omega$ is nonexistent, and must not participate in the approximation of any function on $\Omega$.
- Lastly, when $x$ locates outside $\Omega$ ($x \in \Lambda$), i.e., $\chi(x) = 0$, all the shape functions of $x_I$ and their derivatives evaluated at such $x$ are zero: $\Psi_I(x) = \Psi_I^{\nabla_x}(x) = \Psi_I^{\nabla_y}(x) = 0$. If $x_I$ locate inside $\Omega$, this is consistent with the truncated shape functions for nodes near the boundaries. If $x_I$ is also outside $\Omega$, then the shape functions and their derivatives should be nonexistent (like the previous case). The term $[1 - \chi(x)]\mathbb{I}$ in Eq. (48) is added for computational convenience: $M = \mathbb{I}$ outside $\Omega$, only to allow $M^{-1}$ to exist, and therefore $\Psi_I(x) = \Psi_I^{\nabla_x}(x) = \Psi_I^{\nabla_y}(x) = 0$ is obtained; otherwise, the shape functions on $\Lambda$ become undefined and problematic in computations.

Given the modified forms, we can recover the RKPM discretization for non-periodic problems on $\Omega$, while solving on the extended periodic domain $\mathbb{T}$. This allows us to exploit the RKPM convolutional structure (see the following) via efficient FFT algorithms and the fast circular convolution given by Eq. (34).

### 4.2. Convolutional forms

In the following, we discuss the convolutional forms for the moment matrix $M$, the internal force $\mathbf{f}^{\mathrm{int}}(\mathbf{d}) = \mathbf{K}\mathbf{d}$, the external force $\mathbf{f}^{\mathrm{r}}(\mathbf{r})$, and the solution $\mathbf{u}^h(\mathbf{d})$. Convolutional structures for other terms associated with boundary conditions, nonlinear problems and time dependent problems are briefly discussed in Appendix A, Appendix B, and Appendix C, respectively.

#### 4.2.1. The Moment Matrix

To arrive at the convolutional forms, it is often more convenient to work with indicial notation of vectors and matrices. Using indicial notation for Eq. (48) one gets:

$$\begin{aligned} M_{pq}(x_I) &= [1 - \chi(x_I)]\delta_{pq} + \chi(x_I) \sum_{I=1}^{N} \chi(x_J) H_p(x_I - x_J) H_q^a(x_I - x_J) \\ &= (1 - \chi_I)\delta_{pq} + \chi_I \sum_{I=1}^{N} \chi_J (H_p H_q^a)_{I-J} \\ &= (1 - \chi_I)\delta_{pq} + \chi_I [\chi \circledast H_p H_q^a]_I \quad , \text{and } p, q = 1, 2, \ldots, s \end{aligned} \quad (49)$$

where $\delta_{pq}$ denotes the Kronecker delta. As seen in Eq. (49), subscripts $(I)$ and $(I-J)$ are used to denote the arguments $(x_I)$ and $(x_I - x_J)$.

### 4.2.2. The internal force

This is the most expensive term to be computed within the traditional RKPM. In FC-RKPM, we do not compute the stiffness matrix **K** a priori, and then compute the **Kd** product; instead, a direct approach is adopted to compute $\mathbf{f}^{int} = \mathbf{Kd}$ as a single term.

First, we expand $\mathbf{f}^{int}$ using Eqs. (21), (24), (26), and the Implicit Gradient approximation of the shape function derivatives in Eq. (14):

$$\mathbf{f}_I^{int} = \sum_{J=1}^{N} K_{IJ} d_J = \sum_{J=1}^{N} \left( \int_{\mathbb{T}} \left( \Psi_I^{\nabla_x}(x) \Psi_J^{\nabla_x}(x) + \Psi_I^{\nabla_y}(x) \Psi_J^{\nabla_y}(x) \right) d\mathbb{T} \right) d_J \tag{50}$$

$$= \int_{\mathbb{T}} \left[ \sum_{J=1}^{N} \left( \Psi_I^{\nabla_x}(x) \Psi_J^{\nabla_x}(x) + \Psi_I^{\nabla_y}(x) \Psi_J^{\nabla_y}(x) \right) d_J \right] d\mathbb{T}$$

$$= \int_{\mathbb{T}} \left[ \sum_{J=1}^{N} \Psi_I^{\nabla_x}(x) \Psi_J^{\nabla_x}(x) d_J \right] d\mathbb{T} + \int_{\mathbb{T}} \left[ \sum_{J=1}^{N} \Psi_I^{\nabla_y}(x) \Psi_J^{\nabla_y}(x) d_J \right] d\mathbb{T}$$

Substituting Eqs. (42) and (43):

$$\mathbf{f}_I^{int} = \int_{\mathbb{T}} \left[ \sum_J \chi(x)\chi(x_I)[b^x]^T(x) H^a(x - x_I) \chi(x)\chi(x_J)[b^x]^T(x) H^a(x - x_J) d_J \right] d\mathbb{T} \tag{51}$$

$$+ \int_{\mathbb{T}} \left[ \sum_J \chi(x)\chi(x_I)[b^y]^T(x) H^a(x - x_I) \chi(x)\chi(x_J)[b^y]^T(x) H^a(x - x_J) d_J \right] d\mathbb{T}$$

Taking functions of $x$ and $x_I$ out of the summation $\Sigma_J$ and using indicial notations for vectors, and given that $\chi^2(x) = \chi(x)$ one gets:

$$\mathbf{f}_I^{int} = \int_{\mathbb{T}} \left[ \chi(x)\chi(x_I) b_p^x(x) H_p^a(x - x_I) b_q^x(x) \sum_J H_q^a(x - x_J) \chi(x_J) d_J \right] d\mathbb{T} \tag{52}$$

$$+ \int_{\mathbb{T}} \left[ \chi(x)\chi(x_I) b_p^y(x) H_p^a(x - x_I) b_q^y(x) \sum_J H_q^a(x - x_J) \chi(x_J) d_J \right] d\mathbb{T}$$

, and $p, q = 1, 2, \ldots, s$

where $\sum_J H^a(x - x_J) \chi_J d_J$ is a circular convolution sum on $\mathbb{T}$:

$$\mathbf{f}_I^{int} = \int_{\mathbb{T}} [\chi(x)\chi(x_I) b_p^x(x) H_p^a(x - x_I) b_q^x(x) [\chi d \circledast H_q^a](x)] d\mathbb{T} \tag{53}$$

$$+ \int_{\mathbb{T}} [\chi(x)\chi(x_I) b_p^y(x) H_p^a(x - x_I) b_q^y(x) [\chi d \circledast H_q^a](x)] d\mathbb{T}$$

We use *Direct Nodal Integration* (DNI) [19] to discretize the integrals in Eq. (53). In DNI, quadrature points are taken to be the same as nodes: $x_I$.

$$f_I^{\text{int}} \cong \sum_{S=1}^{N} \chi(x_S)\chi(x_I)b_p^x(x_S)H_p^a(x_S - x_I)b_q^x(x_S)[\chi d \circledast H_q^a](x_S)V(x_S) \qquad (54)$$

$$+ \sum_{S=1}^{N} \chi(x_S)\chi(x_I)b_p^y(x_S)H_p^a(x_S - x_I)b_q^y(x_S)[\chi d \circledast H_q^a](x_S)V(x_S)$$

Where $V(x_S)$ is the quadrature weight, i.e., the volume associated with the node $x_S$. Let:

$$\overline{H^a}(x) = H^a(-x). \qquad (55)$$

Eq. (54) yields:

$$f_I^{\text{int}} = \sum_{S=1}^{N} \chi_S \chi_I (b_p^x)_S (b_q^x)_S (\chi d \circledast H_q^a)_S V_S (\overline{H_p^a})_{I-S} \qquad (56)$$

$$+ \sum_{S=1}^{N} \chi_S \chi_I (b_p^y)_S (b_q^y)_S (\chi d \circledast H_q^a)_S V_S (\overline{H_p^a})_{I-S}$$

Factoring out $\chi_I$ from the summations, we get the final convolutional form for $\mathbf{f}^{\text{int}} = \mathbf{Kd}$ as:

$$f_I^{\text{int}} = \chi_I \left( \{[\chi b_q^x (\chi d \circledast H_q^a) b_p^x V] \circledast \overline{H_p^a}\}_I + \{[\chi b_q^y (\chi d \circledast H_q^a) b_p^y V] \circledast \overline{H_p^a}\}_I \right) \qquad (57)$$

, and $p, q = 1, 2, \ldots, s$.
Note that $p$ and $q$ in Eq. (57) are dummy indices, and therefore, denote the summation $\sum_{q=1}^{s} \sum_{p=1}^{s}$.

### 4.2.3. The external force

Substituting Eq. (23) in Eq. (27), and using Eq. (41) yields:

$$f_I^r = \int_{\mathbb{T}} \Psi_I(x) r(x) \, d\mathbb{T} = \int_{\mathbb{T}} \chi(x)\chi(x_I)[b^0]^{\text{T}}(x) H^a(x - x_I) r(x) \, d\mathbb{T} \qquad (58)$$

Using DNI for quadrature and factoring out $\chi_I$:

$$f_I^r = \chi_I \sum_{S=1}^{N} \chi_S [b^0]_S^{\text{T}} H_{S-I}^a r_S V_S \qquad (59)$$

Using Eq. (55), and indicial notation for vector product, the convolutional form is obtained:

$$f_I^r = \chi_I \sum_{S=1}^{N} \chi_S r_S V_S (b_p^0)_S (\overline{H_p^a})_{I-S} = \chi_I (\chi r V b_p^0 \circledast \overline{H_p^a})_I \qquad , \text{and } p, q = 1, 2, \ldots, s \qquad (60)$$

### 4.2.4. Evaluation of the approximated functions

In most cases, one needs to evaluate the approximated field variable $u^h(x)$, from the coefficients $d_I$ and Eq. (1). Fast convolution sum is applicable to efficient evaluation of $u^h(x)$ as well:

$$u^h(x) = \sum_{J=1}^{N} \Psi_J(x) d_J = \sum_{J=1}^{N} \chi(x)\chi(x_J)[b^0]^{\text{T}}(x) H^a(x - x_J) d_J \qquad (61)$$

The discrete convolution form is easily achieved for the field variable at the nodes:

$$u^h(x_I) = \sum_{J=1}^{N} \Psi_J(x_I)d_J = \sum_{J=1}^{N} \chi(x_I)\chi(x_J)[\boldsymbol{b}^0]^T(x_I)\boldsymbol{H}^a(x_I - x_J)d_J \qquad (62)$$

$$= \chi_I[\boldsymbol{b}^0]_I^T \sum_{J=1}^{N} \chi_J d_J \boldsymbol{H}^a_{I-J}$$

Using indicial notations for vectors:

$$u_I^h = \chi_I(b_p^0)_I \sum_{J=1}^{N} \chi_J d_J (H_p^a)_{I-J} = \chi_I(b_p^0)_I(\chi d \circledast H_p^a)_I \quad , \text{and } p = 1,2,\dots,s \qquad (63)$$

*Remark:* If the interpolated field variable is needed on a position other than the nodes, i.e., on $x \neq x_I$ (for all $I$), the fast convolution using uniform FFT cannot be used. Using non-uniform FFT instead, may be applicable to perform the fast convolution, but this has not been thoroughly investigated, and remains a possibility for future research.

*Remark:* For periodic problems as a special case, one can set $\chi(x) = 1$ for all $x$.

*Remark:* In this study, in order to reach the convolutional forms, DNI has been used for quadrature. Meshfree nodes and quadrature points being identical leads to convolution sums of data collected on identical spatial coordinates. If spatial coordinates of data differ for the two convolving functions, the conventional definition of a convolution sum (regardless of the uniformity of nodal spacing) does not hold. Exploiting fast convolution sums for quadrature methods other than DNI remains a topic for future research.

*Remark*: In certain problems, DNI leads to instabilities that originate in excitation of zero energy modes [19] and non-convergence of the solution [57, 58]. To remedy this issue, numerous treatments have been proposed such as [11, 15, 47]. The convolutional form for RKPM with stabilizers is also left for future extension of the method, but the techniques here can easily be extended to purely node-based natural stabilization [11] and node-based corrections for accuracy [58].

*Remark:* The scalar Poisson problem (and diffusion problem in the Appendix C) are used here as simple model problems for demonstration of the new method. The introduced FC-RKPM is however general, and can be applied to the RKPM discretization of any PDEs, such as equation of motion, electrostatics, mass transfer, etc.; the principles are the same. While, some discussion are provided to a variety of terms arising in the weak form of different PDEs in the appendices, detailed implementation of the method for any of these set of problems can be a topic of future research.

*Remark:* The *fast-convolving kernel* approach is general and applicable to other meshfree methods that employ kernel-type approximations (e.g., smoothed particle hydrodynamics, reproducing kernel collocation method, etc.). The method has been already used for peridynamics [34–37]. At some point in discretization, these meshfree methods need evaluating convolution sums. The fast convolution method is then applicable if the domain is extended to a box, and appropriate corrections via characteristic functions are considered such that the box is partitioned to specific sub-domains.

## 5. Implementation

In this section, we provide the detailed implementation procedure for using FC-RKPM. We discuss: the specifics of extending the domain at the discrete level, necessary adjustments to certain convolving functions, the data structure for programming, and suggested algorithms.

### 5.1. Domain extension

If the domain $\Omega$ is not rectangular, first, one needs to define an enclosing box within which $\Omega$ fits. The Cartesian coordinate system should be aligned with the box edges. It is best to choose a box/coordinate direction that minimizes the "gap", and consequently leads to less total degrees of freedom after discretization considering the same grid spacing. For example, Figure 3 illustrates a good and a poor choice of coordinate directions for a specific $\Omega$. A rule of thumb is to have one coordinate direction aligned with the longest dimension of the geometry.

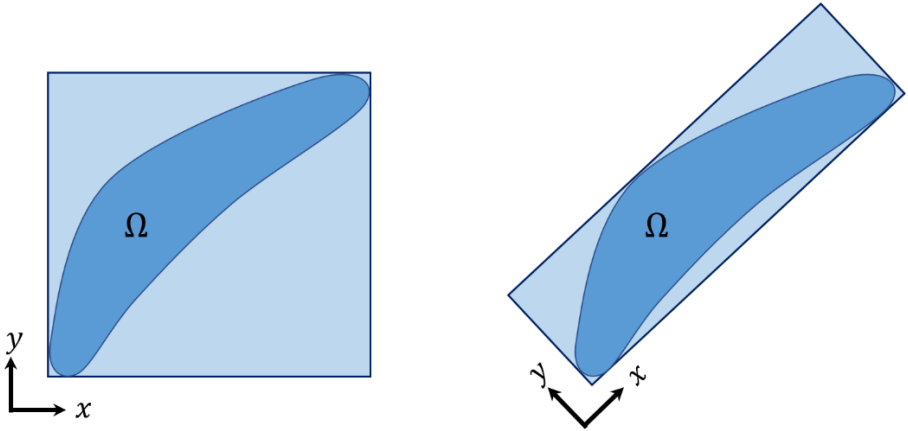

Figure 3. Schematic demonstration of a good (right) and a poor (left) choice of coordinate system for defining the box for a non-rectangular domain $\Omega$.

Given a domain of interest $\Omega$, and a coordinate system, the fitted box is extended in one direction of each coordinate axis (either positive or negative direction). For example, see the extensions $l_{e_x}$ and $l_{e_y}$ for a 2D domain in Figure 4.

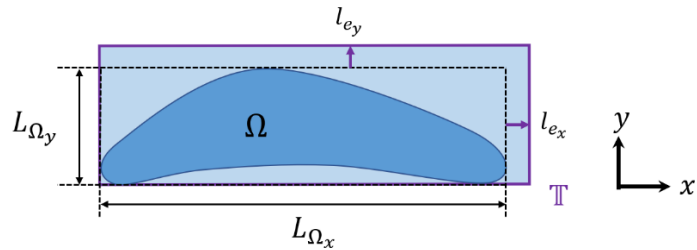

Figure 4. Extension of a 2D generic $\Omega$ by $l_{e_x}$ and $l_{e_y}$ in the two directions.

The minimum acceptable extension for FC-RKPM is to have $l_{e_x}$ and $l_{e_y}$ larger than the support sizes, $a_x$ and $a_y$ respectively. The support sizes are defined based on the discretization:

$$a_x = \widetilde{a_x}\Delta x; \quad a_y = \widetilde{a_y}\Delta y, \tag{64}$$

where $\widetilde{a_x}$ and $\widetilde{a_y}$ are constants. Let:

$$m_x = \lfloor \widetilde{a_x} \rfloor; \quad m_y = \lfloor \widetilde{a_y} \rfloor \qquad (65)$$

where $\lfloor \cdot \rfloor$ denotes flooring (rounding down) operator. Then one gets:

$$l_{e_x} = (m_x + 1)\Delta x \text{ and } l_{e_y} = (m_y + 1)\Delta y \qquad (66)$$

Given the nodal spacing or discretization of $\Omega$, one can use Eq. (66) to find the extensions and define the periodic box. Note that, FFT shows its optimum performance if $N_x$ and $N_y$ (nodes in the whole box, not just $\Omega$) are powers of 2. It is, then, recommended to choose $N_x$ and $N_y$ first, and find the nodal spacing, support sizes, and extensions accordingly:

Let $L_{\Omega_x}$ and $L_{\Omega_y}$ be the largest dimension of $\Omega$ along $x$ and $y$ axes, i.e., the dimension of the fitted box before extension in Figure *4*. Given $N_x = 2^{P_x}$ and $N_y = 2^{P_y}$ ($P_x$ and $P_y$ are positive integers), and substituting $\Delta x = \frac{(L_{\Omega_x} + l_{e_x})}{N_x}$ and $\Delta y = \frac{(L_{\Omega_y} + l_{e_y})}{N_y}$ in Eq. (66), one gets:

$$l_{e_x} = \frac{m_x + 1}{N_x - m_x - 1} L_{\Omega_x} \text{ and } l_{e_y} = \frac{m_y + 1}{N_y - m_y - 1} L_{\Omega_y} \qquad (67)$$

*Remark*: If $\Delta x$ and $\Delta y$ are determined in advance, one can use Eqs. (*64*) to (*66*) to find the extensions. The total number of nodes, however, may not end up as a power of 2. This could slow down the FFT operation a few times, but the method would still be significantly more efficient compared with the traditional approach. This statement is confirmed by results in Section 6.2.2.

### 5.2. Spatial adjustment of the kernel and the monomial basis

For the fast convolution in Eq. (*34*) to be equivalent to the circular convolution sum in Eq. (*39*) on a periodic domain $\mathbb{T}$, it is necessary to have the convolving functions positioned in accordance with the coordinate in which the domain $\mathbb{T}$ is defined. In 2D for example, assume at $\mathbb{T} = [x_{\min} \ x_{\max}] \times [y_{\min} \ y_{\max}]$ where $\boldsymbol{x}_{\min} = \begin{Bmatrix} x_{\min} \\ y_{\min} \end{Bmatrix} \neq \begin{Bmatrix} 0 \\ 0 \end{Bmatrix}$. The adjustment requires to ensure periodicity of all convolving functions, and also shifting the functions such that their zeros coincide with the coordinate $\begin{Bmatrix} x_{\min} \\ y_{\min} \end{Bmatrix}$. To visually demonstrate the procedure, the adjustment for a generic radial kernel function $\phi$, is schematically shown in Figure 5 (1D) and Figure 6 (2D).

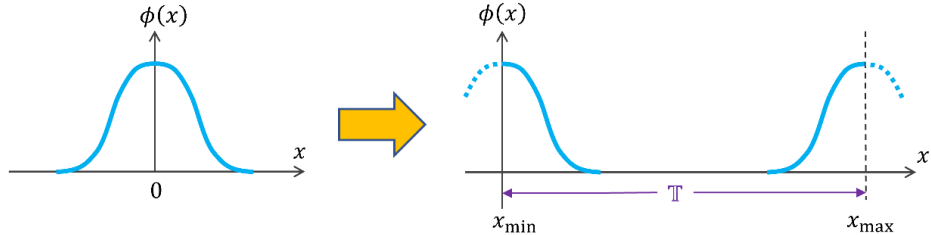

Figure 5. Schematic demonstration of adjusting a kernel function for fast convolution in 1D. Left: kernel function; right: adjusted kernel function on $\mathbb{T}$.

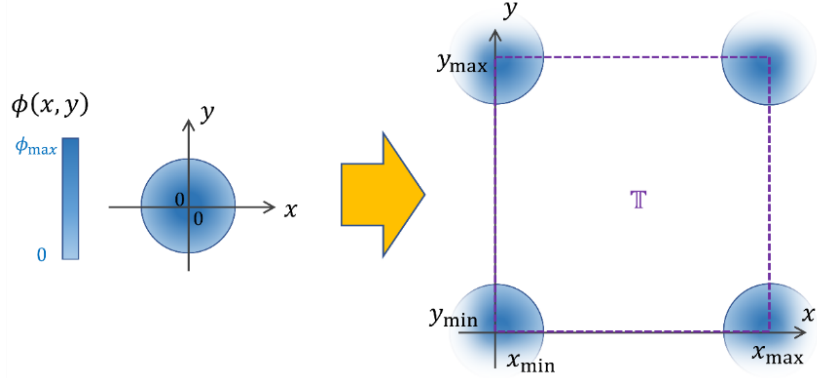

Figure 6. Schematic demonstration of a kernel function (left) and its adjusted version on $\mathbb{T}$ (right) for fast circular convolution in 2D.

One way to create the adjusted functions is to split the functions on the finite size domain centered at zero (with the same size as one period of $\mathbb{T}$), and reorder and shift the split parts, such that the adjusted form on $\mathbb{T}$ is obtained. Figure 7 schematically shows this procedure for a generic function in 2D.

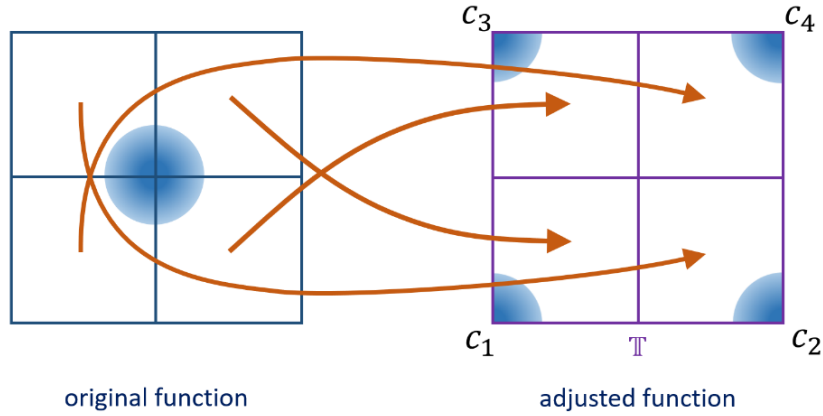

original function      adjusted function

Figure 7. Schematic demonstration of adjusting a generic function for fast circular convolution in 2D. Left: original form of the function; right: adjusted function on $\mathbb{T}$.

The mathematical expression of this adjustment is:

$$\phi^{(adj)}(x) = \sum_{i=1}^{n} \chi_{C_i}(x)\, \phi(x - x_{c_i}); \ \ n = 2^d \tag{68}$$

where $d$ is the spatial dimension, $c_i$ refers to the $i^{th}$ "corner" of the periodic box $\mathbb{T}$ (4 corners in 2D), and $\chi_{C_i}$ is a characteristic function that represents the one $n^{th}$ partition of $\mathbb{T}$ associated with the corner $c_i$:

$$\chi_{C_i}(x) = \begin{cases} 1 & x \in \text{partition } i \\ 0 & \text{else} \end{cases} \tag{69}$$

In the convolutional forms presented in the previous section, functions that need adjustment are: $H(x)$, $H^a(x)$ and $\overline{H^a}(x)$. The latter two can either be directly adjusted or that they can be obtained from Eqs. (44) and (55) after adjustment of $\phi_a$ and $H$.

### 5.3. Data structure and fast convolutions

In this part, the implementation form of the convolutions derived in Section 4.2, and the associated arrays and data structure are presented. Similar to the rest of this article, we focus on the 2D case for demonstration. The 3D case is easily followed by considering one extra dimension.

Let $\mathbb{T} = [x_{\min}\ x_{\max}] \times [y_{\min}\ y_{\max}]$ denote a periodic box, with the uniform discretization given below:

$$\boldsymbol{x}_{ij} = \begin{Bmatrix} x_i \\ y_j \end{Bmatrix} = \begin{Bmatrix} x_{\min} + (i-1)\Delta x \\ y_{\min} + (j-1)\Delta y \end{Bmatrix}; \Delta x = \frac{L_1}{N_1}, \Delta y = \frac{L_2}{N_2}; \tag{70}$$

and $i = 1,2,\ldots,N_1 \quad j = 1,2,\ldots,N_2$, $\quad L_1 = x_{\max} - x_{\min}, \quad L_2 = y_{\max} - y_{\min}$

Let the $x$ and the $y$ coordinates of all nodes, be stored is two distinct 2D arrays of size $N_1 \times N_2$:

$$\mathbf{X} = [X_{ij}]_{N_1 \times N_2}, \text{ and } \mathbf{Y} = [Y_{ij}]_{N_1 \times N_2}; \ i = 1,2,\ldots,N_1 \quad j = 1,2,\ldots,N_2 \tag{71}$$

where $X_{ij}$ denotes the $x$-coordinate of node $\boldsymbol{x}_{ij}$, while $Y_{ij}$ denotes its $y$-coordinate. Here we use the dimensional indexing, because we need arrays to have dimensions consistent with the spatial dimension of the problem. The boldface nonitalic letters denote the $N_1 \times N_2$ 2D implementation arrays.

The discrete version of the monomial basis functions is defined as:

$$\mathbf{H}_p = \mathbb{H}_p(\mathbf{X}, \mathbf{Y}); \text{ and } p = 1,2,\ldots,s \tag{72}$$

where the function $\mathbb{H}_p$ takes the nodal coordinate arrays $\mathbf{X}$ and $\mathbf{Y}$ as inputs, and returns an $N_1 \times N_2$ array containing the $p$-th element of the vector $\boldsymbol{H}(\boldsymbol{x})$ in Eq. (4), evaluated for all nodes:

$$\mathbf{H}_p = \left[(H_p)_{ij}\right]_{N_1 \times N_2}; \ i = 1,2,\ldots,N_1 \quad j = 1,2,\ldots,N_2 \tag{73}$$

and $(H_p)_{ij} = H_p(\boldsymbol{x}_{ij})$. For example, if $\boldsymbol{H}(\boldsymbol{x}_{ij}) = [1, x_i, y_j]^\mathrm{T}$ is the vector of basis functions, the corresponding $\mathbf{H}_p$ arrays are:

$$\mathbf{H}_1 = [(H_1)_{ij}] = \mathbf{1} \quad (\forall i,j: (H_1)_{ij} = 1) \tag{74}$$

$$\mathbf{H}_2 = [(H_2)_{ij}] = \mathbf{X} \text{ (see Eq. (71))}$$

$$\mathbf{H}_3 = [(H_3)_{ij}] = \mathbf{Y} \text{ (see Eq. (71))}$$

The discrete form of the kernel function is defined as well:

$$\boldsymbol{\Phi}_a = \Phi(\mathbf{X}, \mathbf{Y}) \tag{75}$$

where the discrete function $\Phi$ takes the nodal coordinates arrays $\mathbf{X}$ and $\mathbf{Y}$ as input, and returns an $N_1 \times N_2$ array containing the nodal kernel values:

$$\boldsymbol{\Phi}_a = \left[(\phi_a)_{ij}\right]_{N_1 \times N_2} \tag{76}$$

$(\phi_a)_{ij} = \phi_a(\boldsymbol{x}_{ij})$ is the kernel function evaluated at the node $\boldsymbol{x}_{ij}$.

Given the definitions for $H_p^a$ and $\overline{H_p^a}$ in Eqs. (44) and (55), the discrete functions above, and the adjustment procedure described in Eq. (68), the $N_1 \times N_2$ arrays for the adjusted $H_p^a$ and $\overline{H_p^a}$ functions are created as follows (needs to be done for each $p = 1,2,\ldots,s$):

$$\mathbf{H}_p^a = \sum_{i=1}^{4} \boldsymbol{\chi}_{C_i}(\mathbf{X},\mathbf{Y}) \circ \mathbb{H}_p(\mathbf{X} - x_{c_i}, \mathbf{Y} - y_{c_i}) \circ \Phi(\mathbf{X} - x_{c_i}, \mathbf{Y} - y_{c_i}) \tag{77}$$

$$\overline{\mathbf{H}_p^a} = \sum_{i=1}^{4} \boldsymbol{\chi}_{C_i}(\mathbf{X},\mathbf{Y}) \circ \mathbb{H}_p\left(-(\mathbf{X} - x_{c_i}), -(\mathbf{Y} - y_{c_i})\right) \circ \Phi\left(-(\mathbf{X} - x_{c_i}), -(\mathbf{Y} - y_{c_i})\right) \tag{78}$$

Where $\boldsymbol{\chi}_{C_i}(\mathbf{X},\mathbf{Y})$ is the discrete version of the characteristic function defined by Eq. (69), which is an $N_1 \times N_2$ array containing 1 for $\boldsymbol{x}_{ij}$ that locate on the partition corresponding to $c_i$ and 0 otherwise. The operator ($\circ$) in equations above denotes the elementwise product of two arrays of the same dimension.

The arrays $\mathbf{H}_p$ ($p = 1,2,\ldots,s$) need to be adjusted as well:

$$\mathbf{H}_p = \sum_{i=1}^{4} \boldsymbol{\chi}_{C_i}(\mathbf{X},\mathbf{Y}) \circ \mathbb{H}_p(\mathbf{X} - x_{c_i}, \mathbf{Y} - y_{c_i}) \tag{79}$$

Having the necessary arrays computed and stored, the implementation for fast evaluation of the terms derived in Section 4.2 is presented as follows:

*Moment matrix and its inverse*

Using the fast convolution sum given by Eq. (*34*), Eq. (49) can be implemented as:

$$\mathbf{M}_{pq} = (1 - \boldsymbol{\chi})\delta_{pq} + \boldsymbol{\chi} \circ \mathbf{F}^{-1}\left[\mathbf{F}(\boldsymbol{\chi}) \circ \mathbf{F}(\mathbf{H}_p \circ \mathbf{H}_q^a)\right] \quad , \text{and } p,q = 1,2,\ldots,s \tag{80}$$

Where $\mathbf{F}$ and $\mathbf{F}^{-1}$ denote the FFT and inverse FFT operations, and:

$$\boldsymbol{\chi} = [\chi_{ij}]_{N_1 \times N_2} = [\chi(\boldsymbol{x}_{ij})]_{N_1 \times N_2} \tag{81}$$

$$\mathbf{M}_{pq} = \left[(M_{pq})_{ij}\right]_{N_1 \times N_2} = [M_{pq}(\boldsymbol{x}_{ij})]_{N_1 \times N_2} \tag{82}$$

are $N_1 \times N_2$ arrays ($i = 1,2,\ldots,N_1$; $j = 1,2,\ldots,N_2$). $\mathbf{H}_p$ and $\mathbf{H}_q^a$ are given by Eqs. (79) and (77). Note that similar to $\mathbf{H}_p$, $\mathbf{M}_{pq}$ is a distinct 2D array for each $p$ and $q$. For example, if $s = 3$, $\mathbf{M}_{11}, \mathbf{M}_{12}, \ldots, \mathbf{M}_{33}$ are each a 2D $N_1 \times N_2$ array.

Once $\mathbf{M}_{pq}$ are found from Eq. (*80*), the $s \times s$ moment matrix $\mathbf{M}(\boldsymbol{x}_{ij})$ is assembled for each node $\boldsymbol{x}_{ij}$, and its inverse is computed: $\mathbf{M}^{-1}(\boldsymbol{x}_{ij})$. This is shown in Algorithm 1 in the next section.

Having $\mathbf{M}^{-1}(\boldsymbol{x}_{ij})$, the following $N_1 \times N_2$ arrays are constructed to be used in the FC-RKPM computations:

$$\mathbf{b}_p^0 = \left[(b_p^0)_{ij}\right]_{N_1 \times N_2} = [b_p^0(\boldsymbol{x}_{ij})]_{N_1 \times N_2}; \text{ and } b_p^0(\boldsymbol{x}_{ij}) = [M^{-1}]_{1p}(\boldsymbol{x}_{ij}), \tag{83}$$

$$\mathbf{b}_p^x = \left[(b_p^x)_{ij}\right]_{N_1 \times N_2} = [b_p^x(\boldsymbol{x}_{ij})]_{N_1 \times N_2}; \text{ and } b_p^x(\boldsymbol{x}_{ij}) = -[M^{-1}]_{2p}(\boldsymbol{x}_{ij}), \tag{84}$$

$$\mathbf{b}_p^y = \left[(b_p^y)_{ij}\right]_{N_1 \times N_2} = [b_p^y(\boldsymbol{x}_{ij})]_{N_1 \times N_2}; \text{ and } b_p^y(x_i, y_j) = -[M^{-1}]_{3p}(\boldsymbol{x}_{ij}), \tag{85}$$

for each $p = 1,2,...,s$.

*Internal force*

Using the fast convolution sum in Eq. (*34*), Eq. (57) can be implemented as:

$$\mathbf{f}^{\text{int}} = \boldsymbol{\chi} \circ \{\mathbf{F}^{-1}[\mathbf{F}(\boldsymbol{\chi} \circ \mathbf{V} \circ \mathbf{b}_p^x \circ \mathbf{b}_q^x \circ \mathbf{F}^{-1}[\mathbf{F}(\boldsymbol{\chi} \circ \mathbf{d}) \circ \mathbf{F}(\mathbf{H}_q^a)]) \circ \mathbf{F}(\overline{\mathbf{H}_p^a})] \qquad (86)$$
$$+ \mathbf{F}^{-1}[\mathbf{F}(\boldsymbol{\chi} \circ \mathbf{V} \circ \mathbf{b}_p^y \circ \mathbf{b}_q^y \circ \mathbf{F}^{-1}[\mathbf{F}(\boldsymbol{\chi} \circ \mathbf{d}) \circ \mathbf{F}(\mathbf{H}_q^a)]) \circ \mathbf{F}(\overline{\mathbf{H}_p^a})]\}$$

, and $p, q = 1,2,...,s$.
where:

$$\mathbf{V} = [V_{ij}]_{N_1 \times N_2} = [V(\boldsymbol{x}_{ij})]_{N_1 \times N_2} \qquad (87)$$

$$\mathbf{d} = [d_{ij}]_{N_1 \times N_2} = [d(\boldsymbol{x}_{ij})]_{N_1 \times N_2} \qquad (88)$$

$$\overline{\mathbf{H}_p^a} = [(\overline{H_p^a})_{ij}]_{N_1 \times N_2} = [\overline{H_p^a}(\boldsymbol{x}_{ij})]_{N_1 \times N_2} \quad , \text{for } p = 1,2,...,s \qquad (89)$$

are $N_1 \times N_2$ arrays ($i = 1,2,...,N_1 \quad j = 1,2,...,N_2$).

Note that for nodes that locate on the interior of $\Omega$, $V_{ij} = \Delta x \Delta y$ for uniform nodal spacing (fixed $\Delta x$ and $\Delta y$). For nodes on the boundary $V_{ij}$ is a fraction of $\Delta x \Delta y$ depending on the boundary's shape on that node.

Eq. (86) can be further simplified. Defining:

$$\mathbf{A}^x = \mathbf{b}_q^x \circ \mathbf{F}^{-1}[\mathbf{F}(\boldsymbol{\chi} \circ \mathbf{d}) \circ \mathbf{F}(\mathbf{H}_q^a)] \qquad (90)$$
$$\mathbf{A}^y = \mathbf{b}_q^y \circ \mathbf{F}^{-1}[\mathbf{F}(\boldsymbol{\chi} \circ \mathbf{d}) \circ \mathbf{F}(\mathbf{H}_q^a)]$$

(with $q = 1,2,...,s$), and factoring out $\mathbf{F}(\overline{\mathbf{H}_p^a})$ from Eq. (*86*), one then gets:

$$\mathbf{f}^{\text{int}} = \boldsymbol{\chi} \circ \mathbf{F}^{-1}\{[\mathbf{F}(\boldsymbol{\chi} \circ \mathbf{V} \circ \mathbf{b}_p^x \circ \mathbf{A}^x) + \mathbf{F}(\boldsymbol{\chi} \circ \mathbf{V} \circ \mathbf{b}_p^y \circ \mathbf{A}^y)] \circ \mathbf{F}(\overline{\mathbf{H}_p^a})\}, \quad \text{and } p = 1,2,...,s. \qquad (91)$$

*External force*

Eq. (*60*) is implemented as:

$$\mathbf{f}^r = \boldsymbol{\chi} \circ \mathbf{F}^{-1}[\mathbf{F}(\boldsymbol{\chi} \circ \mathbf{V} \circ \mathbf{b}_p^0 \circ \mathbf{r}) \circ \mathbf{F}(\overline{\mathbf{H}_p^a})] \qquad (92)$$
, and $p = 1,2,...,s$.
Where

$$\mathbf{r} = [r_{ij}]_{N_1 \times N_2} = [r(\boldsymbol{x}_{ij})]_{N_1 \times N_2} \qquad (93)$$

*Evaluation of the field variable*

Eq. (63) is implemented as:

$$\mathbf{u}^h = \boldsymbol{\chi} \circ \mathbf{b}_p^0 \circ \mathbf{F}^{-1}[\mathbf{F}(\boldsymbol{\chi} \circ \mathbf{d}) \circ \mathbf{F}(\mathbf{H}_p^a)], \text{and } p = 1,2,...,s. \qquad (94)$$

The computer implementation for fast evaluation of boundary integrals, nonlinear terms, and the terms associated with time-dependent problems are respectively given in Appendix A (Part A.2), Appendix B (Part B.2), and Appendix C (Part C.2).

### 5.4. Solvers and algorithms

In this part, we suggest the general construct for FC-RKPM solvers, and lay out some of the main algorithms.

First, we describe algorithms for particular subroutines that are most likely needed in an FC-RKPM analysis. Similar to the rest of this study, the algorithms are presented for the 2D case; but can be easily modified for 3D or 1D.

#### 5.4.1. Subroutines

Algorithm 1 lays out the subroutine for computing $M(x_{ij})$ and $M^{-1}(x_{ij})$, and storing $\mathbf{b}_p^0, \mathbf{b}_p^x, \mathbf{b}_p^y$ ($p = 1, \ldots, s$) from $M^{-1}$.

---
Algorithm 1.     Subroutine for evaluation of the $M$ and $M^{-1}$, and obtaining $\mathbf{b}_p^0, \mathbf{b}_p^x, \mathbf{b}_p^y$.

Precomputed: $\mathbf{H}_p$; $\mathbf{H}_q^a$

$\hat{\chi} = \mathbf{F}(\chi)$

For $p = 1:s$
    For $q = p:s$ ($M$ is symmetric)
        If $q = p$ then
            $\mathbf{M}_{pq} = (1 - \chi) + \chi \circ \mathbf{F}^{-1}[\hat{\chi} \circ \mathbf{F}(\mathbf{H}_p \circ \mathbf{H}_q^a)]$
        Else
            $\mathbf{M}_{pq} = \chi \circ \mathbf{F}^{-1}[\hat{\chi} \circ \mathbf{F}(\mathbf{H}_p \circ \mathbf{H}_q^a)]$
            $\mathbf{M}_{qp} = \mathbf{M}_{pq}$
        End
    End
End
For $i = 1:N_1$
    For $j = 1:N_2$
        Assemble $M(x_{ij})$ locally: $M(x_{ij}) = \begin{bmatrix} M_{11}(i,j) & \cdots & M_{1q}(i,j) \\ \vdots & \ddots & \vdots \\ M_{p1}(i,j) & \cdots & M_{pq}(i,j) \end{bmatrix}$
        Compute its inverse: $M^{-1}(x_{ij})$
        For $p = 1:s$
            $\mathbf{b}_p^0(i,j) = [M^{-1}]_{1p}(x_{ij})$
            $\mathbf{b}_p^x(i,j) = -[M^{-1}]_{2p}(x_{ij})$
            $\mathbf{b}_p^y(i,j) = -[M^{-1}]_{3p}(x_{ij})$
        End
    End
End

---

This algorithm computes $M$ matrices using $(s^2 + 2)$ FFT/iFFT operations.

Algorithm 2 gives the subroutine for the internal force (linear case) described by Eq. (91).

| Algorithm 2. | Subroutine for evaluation of the linear internal force: $\mathbf{f}^{int}(\mathbf{d})$ |
|---|---|

Precomputed:
$\mathbf{b}_p^x; \mathbf{b}_p^y;\ \mathbf{C}_p^x = \boldsymbol{\chi} \circ \mathbf{V} \circ \mathbf{b}_p^x;\ \mathbf{C}_p^y = \boldsymbol{\chi} \circ \mathbf{V} \circ \mathbf{b}_p^y; \widehat{\mathbf{H}_p^a} = \mathbf{F}(\mathbf{H}_p^a);\ \widehat{\overline{\mathbf{H}_p^a}} = \mathbf{F}(\overline{\mathbf{H}_p^a});$ (for all $p = 1:s$)

Given $\mathbf{d}$ as input:

$\widehat{\mathbf{d}_\chi} = \mathbf{F}(\boldsymbol{\chi} \circ \mathbf{d})$
Set $\mathbf{A}^x = 0;\ \mathbf{A}^y = 0;$
For $p = 1:s$
 $\quad \mathbf{D}_p = \mathbf{F}^{-1}[\widehat{\mathbf{d}_\chi} \circ \widehat{\mathbf{H}_p^a}]$
 $\quad \mathbf{A}^x = \mathbf{A}^x + \mathbf{b}_p^x \circ \mathbf{D}_p$
 $\quad \mathbf{A}^y = \mathbf{A}^y + \mathbf{b}_p^y \circ \mathbf{D}_p$
End

$\widehat{\mathbf{B}} = 0$
For $p = 1:s$
 $\quad \widehat{\mathbf{B}} = \widehat{\mathbf{B}} + \mathbf{F}(\mathbf{C}_p^x \circ \mathbf{A}^x + \mathbf{C}_p^y \circ \mathbf{A}^y) \circ \widehat{\overline{\mathbf{H}_p^a}}$
End

$\mathbf{f}^{int} = \boldsymbol{\chi} \circ \mathbf{F}^{-1}(\widehat{\mathbf{B}})$

With this algorithm, the internal force ($\mathbf{f}^{int}$) for the linear case, regardless of the dimension, is evaluated via $2(s+1)$ FFT/iFFT operations. Given $s = 4$ for a 3D case that uses linear monomial basis function, 10 FFT/iFFT operations are sufficient to compute $\mathbf{f}^{int}$.

For evaluation of $\mathbf{u}^h(\mathbf{d})$ given by Eq. (94) the following subroutine can be used:

| Algorithm 3. | Subroutine for evaluation of the approximated field: $\mathbf{u}^h(\mathbf{d})$ |
|---|---|

Precomputed in initialization: $\mathbf{b}_p^0; \mathbf{b}_p^\chi = \boldsymbol{\chi} \circ \mathbf{b}_p^0; \widehat{\mathbf{H}_p^a} = \mathbf{F}(\mathbf{H}_p^a)$ (for all $p = 1:s$)

Given $\mathbf{d}$

$\widehat{\mathbf{d}_\chi} = \mathbf{F}(\boldsymbol{\chi} \circ \mathbf{d})$
Set $\mathbf{u}^h = 0;$
For $p = 1:s$
 $\quad \mathbf{D}_p = \mathbf{F}^{-1}[\widehat{\mathbf{d}_\chi} \circ \widehat{\mathbf{H}_p^a}]$
 $\quad \mathbf{u}^h = \mathbf{u}^h + \mathbf{b}_p^\chi \circ \mathbf{D}_p$
End

As observed this subroutine requires $(s + 1)$ FFT/iFFT operations.

For the external force given by Eq. (92), the subroutine in Algorithm 4 is provided.

| Algorithm 4. | Subroutine for the external internal force given by Eq. (92): $\mathbf{f}^r(\mathbf{r})$ |
|---|---|

Precomputed in initialization: $\mathbf{b}_p^0$; $\mathbf{C}_p^0 = \boldsymbol{\chi} \circ \mathbf{V} \circ \mathbf{b}_p^0$; $\widehat{\mathbf{H}_p^a} = \mathbf{F}(\mathbf{H}_p^a)$; (for all $p = 1:s$)

Given $\mathbf{r}$ as input:

Set $\widehat{\mathbf{B}} = 0$
For $p = 1:s$
$\quad \widehat{\mathbf{B}} = \widehat{\mathbf{B}} + \mathbf{F}(\mathbf{C}_p^0 \circ \mathbf{r}) \circ \widehat{\mathbf{H}_p^a}$
End

$\mathbf{f}^r = \boldsymbol{\chi} \circ \mathbf{F}^{-1}(\widehat{\mathbf{B}})$

Subroutines for evaluating boundary integrals, nonlinear terms, and terms associated with time dependent problems are respectively provided in Appendix A (Part A.2), Appendix B (Part B.2), and Appendix C (Part C.2).

### 5.4.2. Boundary conditions

As mentioned in Section 2.3 and 2.4, boundary integrals (Naumann and weakly enforced essential BC) can be computed by the traditional approach which is direct quadrature over the specific boundaries. They can also be modified to be computed by FFT-accelerated sums (e.g., see Appendix A (Part A.2) in Appendix A). Regardless of how they are computed, the boundary terms should be structured into $N_1 \times N_2$ arrays using dimensional indexing to be compatible with the rest of the system.

In the case of strong enforcement of essential BC (Eq. (29)) we set $d_{ij} = g(\mathbf{x}_{ij})$ on the boundary nodes ($\mathbf{x}_{ij} \in \Gamma_g$) during the initialization stage (prior to the solver main loop). Then, at each iteration (implicit solvers) or time step (explicit solvers), the coefficients $d_{ij}$ in the array $\mathbf{d}$ are only updated for nodes on $\Omega$; the $d_{ij}$'s on $\Gamma_g$ are left unchanged. The following characteristic functions is helpful for this type of BC enforcement:

$$\chi_{\Gamma_g}(\mathbf{X}, \mathbf{Y}) = \begin{cases} 1 & \mathbf{x}_{ij} \in \Gamma_g \\ 0 & \text{else} \end{cases} \tag{95}$$

$$\chi_\Omega = \chi - \chi_{\Gamma_g} \tag{96}$$

### 5.4.3. Solvers

Given the subroutines for individual terms of the weak forms, the following algorithm is suggested for static/steady state problems. Explicit and implicit solvers for time-dependent problems (e.g., elastodynamics, transient diffusion, etc.) are discussed in Appendix C (Part C.3).

*Static solvers*

In FC-RKPM static analysis *Krylov subspace iterative methods* such as conjugate gradient, GMRES, etc. [59] should be used for solving the linear systems. In each iteration of the Krylov subspace methods, only evaluation of $\mathbf{K}\tilde{\mathbf{d}}$ product is required, where $\mathbf{K}$ is the stiffness matrix of the system and $\tilde{\mathbf{d}}$ is being updated in each iteration. In FC-RKPM, given any $\tilde{\mathbf{d}}$, the corresponding product $\mathbf{K}\tilde{\mathbf{d}}$ can be computed in a wholistic efficient approach as $\tilde{\mathbf{f}}^{int} = \mathbf{f}^{int}(\tilde{\mathbf{d}})$ using Algorithm 2, where the separate construction of $\mathbf{K}$ and the direct matrix vector products are avoided. Direct methods for solving linear systems (e.g., Gaussian elimination) and stationary iterative methods a.k.a. successive relaxation methods (e.g., Jacobi, Gauss-Seidel, etc.) rely on updating the elements of the stiffness matrix $\mathbf{K}$ in each step or iteration, for which the derived convolutional form and the fast operation in Eq. (*91*) are no longer relevant. However, Krylov subspace methods outperform both direct and successive relaxation iterative methods for large systems [60]; they are in fact the best choice for large scale problems, for which the introduced FC-RKPM is primarily intended.

Below is the general outline suggested for time-independent linear analysis:

---
Algorithm 5.    Algorithm outline for linear static (time-independent) analysis.
---

***Read inputs:***
*Discrete data from domain & geometry:* $\mathbf{X}, \mathbf{Y}, \mathbf{V}, \boldsymbol{\chi}, \boldsymbol{\chi}_{\Gamma_g}$;
*Physical parameters:* material properties, boundary conditions' info, body force/heat source $\mathbf{r}(\mathbf{X}, \mathbf{Y})$;
*RKPM parameters:* $s, n, a_x, a_y$;
*Solver parameters:* e.g., tolerance
***Initialize:***
Define functions: $\Phi(\mathbf{X}, \mathbf{Y}), \mathbb{H}_p(\mathbf{X}, \mathbf{Y})$ for $p = 1, \dots, s$
Compute and store: $\mathbf{H}_p, \mathbf{H}_p^a, \overline{\mathbf{H}_p^a}$ (for $p = 1, \dots, s$) from Eqs. (*77*) to (*79*); $\mathbf{b}_p^0, \mathbf{b}_p^x, \mathbf{b}_p^y$ from Algorithm 1; and other arrays needed in the subroutines for computing various terms.
Compute external force (if any): $\mathbf{f}^r$ (Algorithm 4)
Compute boundary terms (if any): e.g., $\mathbf{f}^q$ from Algorithm A.1 (in Appendix A)
Set $\mathbf{d} = 0$ initially
***Solver:*** linear conjugate gradient
While error > tolerance: compute $\mathbf{f}^{int}$ (Algorithm 2) as needed and update $\mathbf{d}$
***Output:***
Compute $\mathbf{u}^h$ (Algorithm 3)
Write outputs

---

For nonlinear static analysis see Algorithm B.3 in Appendix B.

## 6. Numerical examples

In this study, we present numerical examples obtained by FC-RKPM. We use the method of manufactured solution to perform convergence studies in 1D, 2D, and 3D. We also compare the performance of the new method with the traditional RKPM in terms of runtime and memory allocation.

### 6.1. Verification and convergence

Here we study the method's convergence in 1D, 2D, and 3D.

#### 6.1.1. 1D

Consider the following Poisson problem in 1D, subjected to Dirichlet BCs:

$$\begin{cases} \nabla^2 u + 2 = 0, & \text{on } \Omega = [-1,1] \\ u = 0, & \text{on } x = -1,1 \end{cases}, \tag{97}$$

The exact solution to Eq. (97) is:

$$u(x) = 1 - x^2 \tag{98}$$

According to Section 5.1, the domain $\Omega$ needs to be extended with $l_e$ from one end to form $\mathbb{T} = [-1, 1 + l_e)$ as the periodic 1D "box". having $L_\Omega = 2$ and choosing $\tilde{a} = 1.5$, Eqs. (64) to (67) gives:

$$l_e = \frac{4}{N-1} \tag{99}$$

where $N$ is the total number of nodes and selected to be a power of 2.

Once $N$ is selected, arrays **X**, **V**, **χ**, $\boldsymbol{\chi}_{\Gamma_g}$ (size: $N \times 1$) are created as inputs for a MATLAB code that is developed based on the solver in Algorithm 5, which uses a linear conjugate gradient solver to find **d** and $\mathbf{u}^h$. The convergence tolerance for the CG solver is set to $10^{-12}$. Figure 8 shows a visual comparison between the numerical result obtained by FC-RKPM with $N = 2^6$, and the exact solution.

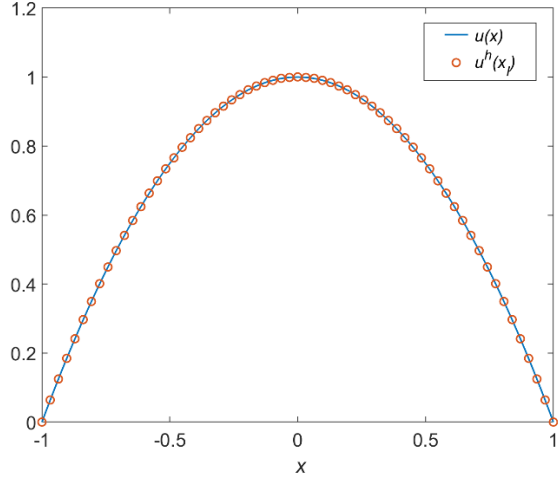

Figure 8. FC-RKPM solution versus exact solution to the 1D Poisson problem.

We performed a convergence study by varying $N = 2^P$ ($P = 3, 6, \ldots, 12$). For the error measure we used the $L_2$-norm error defined by

$$\left| u^h - u \right|_{L_2} = \left( \int_\Omega \left| u^h(x) - u(x) \right|^2 dx \right)^{\frac{1}{2}} \tag{100}$$

We used Gaussian quadrature with 5 Gauss points per cell to compute the norm integral. The convergence study is plotted in Figure 9 and shows a quadratic rate of convergence as one would expect from RKPM using $n = 1$ [61].

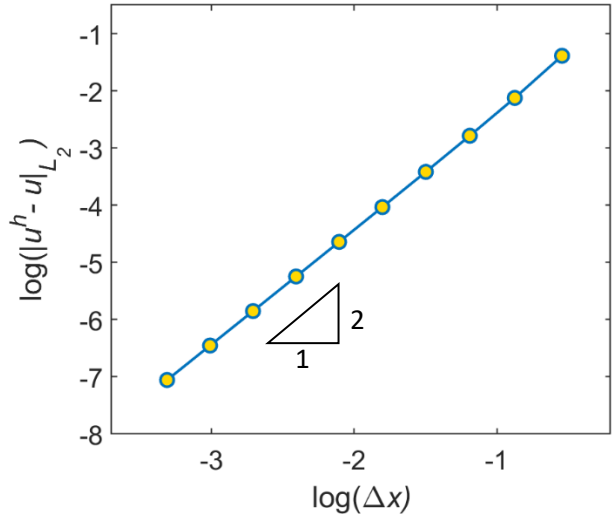

Figure 9. Convergence of FC-RKPM solution for the 1D Poisson problem.

We have also implemented the standard RKPM in a separate MATLAB code which uses the same linear CG for solving the system. The spatial discretization was carried out such that the traditional method uses the same nodes over Ω as those in FC-RKPM. Note that the traditional method does not have the extension. The results from the traditional RKPM were identical to those from the new method, meaning that the errors obtained from both methods were the same up to the machine precision. This shows that new method with fast convolutions and the described modifications, reproduces the solution obtained by the original RKPM.

### 6.1.2. 2D

Consider the following Poisson problem as a 2D example:

$$\begin{cases} \nabla^2 u + 4 - 2x^2 - 2y^2 = 0, & \text{on } \Omega = [-1,1] \times [-1,1] \\ u = 0, & \text{on } x = -1,1; \text{ and on } y = -1,1 \end{cases} \quad (101)$$

with the exact solution being:

$$u(x,y) = (1 - x^2)(1 - y^2) \tag{102}$$

According to Section 5.1, the domain Ω needs to be extended to $\mathbb{T} = \left[-1, 1 + l_{e_x}\right) \times \left[-1, 1 + l_{e_y}\right)$ as the periodic 2D box. having $L_{\Omega_x} = L_{\Omega_y} = 2$ and choosing $\widetilde{a_x} = \widetilde{a_y} = 1.5$, Eqs. (64) to (67) gives:

$$l_{e_x} = \frac{4}{N_x - 2}, \quad \text{and} \quad l_{e_y} = \frac{4}{N_y - 2} \tag{103}$$

where $N_x$ and $N_y$ are selected to be powers of 2. Note that $N_x$ and $N_y$ can be the same or different.

Once $N_x$ and $N_y$ are given, the 2D arrays **X**, **V**, **χ**, $\boldsymbol{\chi}_{\Gamma_g}$ are created as inputs for a 2D MATLAB code following Algorithm 5, using a linear conjugate gradient solver to find **d** and $\mathbf{u}^h$ (Similar to the 1D example). The convergence tolerance for the CG solver is again set to $10^{-12}$. Figure 10 shows a visual comparison between the numerical solution obtained by FC-RKPM with $N_x = N_y = 2^8$ ($N = 2^{16}$), and the exact solution.

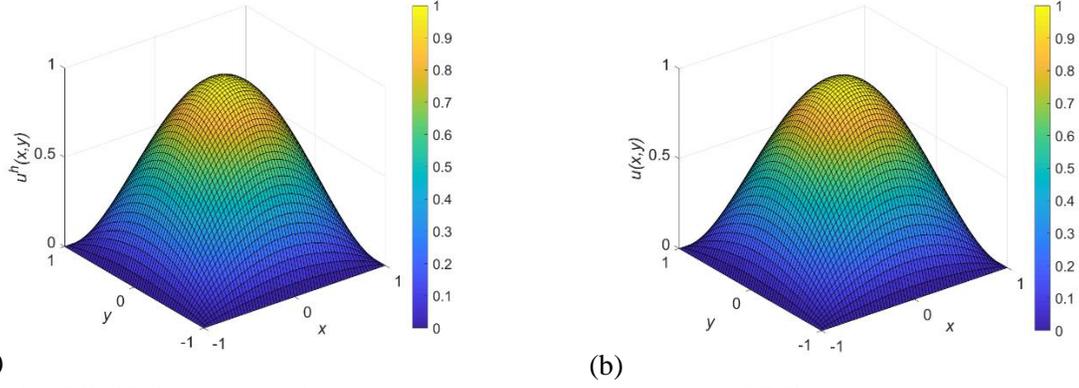

(a)                                                   (b)

Figure 10. FC-RKPM solution (a); and the exact solution (b); for the 2D Poisson problem.

For the convergence study in 2D we varied $N_x$ and $N_y$ from $2^3$ to $2^{10}$ (64 to 1,048,576 nodes), and for the error measures we used the normalized $L_2$-norm and the $L_\infty$-norm of the error at nodes defined by:

$$e_{L_2} = \left( \frac{\sum_{i,j} |u^h(\boldsymbol{x}_{ij}) - u(\boldsymbol{x}_{ij})|^2}{\sum_{i,j} u(\boldsymbol{x}_{ij})^2} \right)^{\frac{1}{2}} \tag{104}$$

$$e_{L_\infty} = \frac{\max_{i,j} |u^h(\boldsymbol{x}_{ij}) - u(\boldsymbol{x}_{ij})|}{\max_{i,j} |u(\boldsymbol{x}_{ij})|} \tag{105}$$

The convergence result is plotted in Figure 11 and confirms the quadratic rate of convergence.

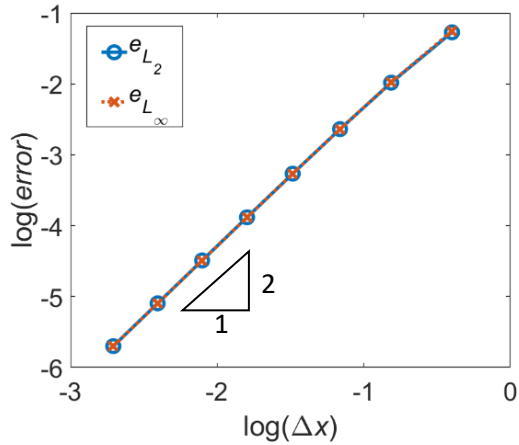

Figure 11. Convergence of the FC-RKPM solution to the exact solution for the 2D Poisson problem.

Again the traditional RKPM with the same nodes on Ω, resulted in identical results as in those in Figure 10 and Figure 11.

$e_{L_2}$ and $e_{L_\infty}$, being close implies that the error is uniformly distributed on Ω and is not localized at a particular point, e.g., near boundaries.

### 6.1.3. 3D

The following Poisson problem is chosen for the 3D example:

$$\begin{cases} \nabla^2 u + 2[3 - 2(x^2 + y^2 + z^2) + x^2y^2 + x^2z^2 + y^2z^2] = 0, & \text{on } \Omega = [-1,1]^3 \\ u = 0, & \text{on } x = -1,1; \text{and on } y = -1,1; \text{and on } z = -1,1 \end{cases} \quad (106)$$

with the exact solution being:

$$u(x,y,z) = (1-x^2)(1-y^2)(1-z^2) \quad (107)$$

The extended domain is then $\mathbb{T} = [-1, 1 + l_{e_x}) \times [-1, 1 + l_{e_y}) \times [-1, 1 + l_{e_z})$ as the periodic 3D box. having $L_{\Omega_x} = L_{\Omega_y} = L_{\Omega_z} = 2$ and choosing $\widetilde{a_x} = \widetilde{a_y} = \widetilde{a_z} = 1.5$, Eqs. (64) to (67) gives:

$$l_{e_x} = \frac{4}{N_x - 2}, \quad l_{e_y} = \frac{4}{N_y - 2}, \quad l_{e_z} = \frac{4}{N_z - 2} \quad (108)$$

where $N_x$, $N_y$, and $N_z$ are selected to be powers of 2.

The 3D problem is solved using a 3D MATLAB code similar to those for the 1D and the 2D examples. Figure 12 shows a visual comparison between the numerical results obtained by FC-RKPM using $N_x = N_y = N_z = 2^6$ ($N = 2^{18}$), and the exact solution.

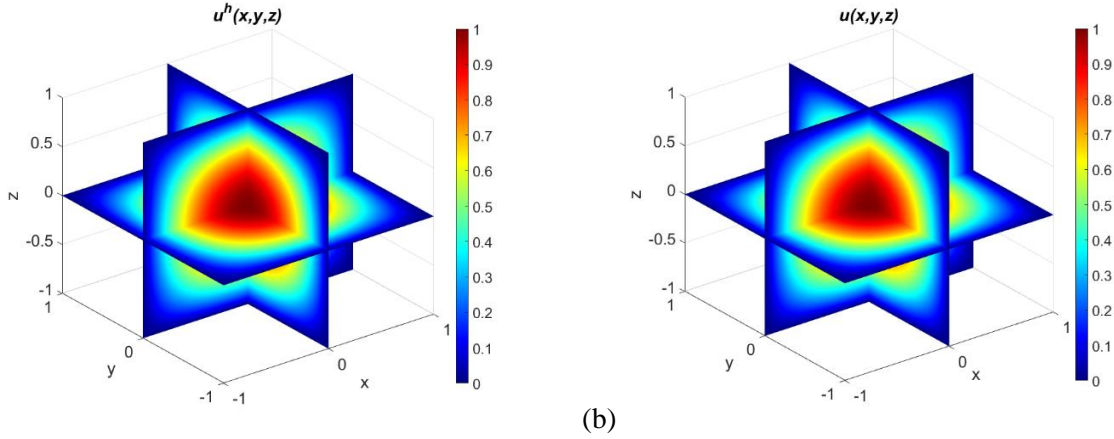

(a)          (b)

Figure 12. FC-RKPM solution (a); and the exact solution (b); for the 2D Poisson problem.

For the convergence study in 3D we changed $N_x$, $N_y$, and $N_z$ from $2^2$ to $2^8$ (64 to 16,777,216 nodes), and for the error measures we used the 3D version of Eqs. (104) and (105). The convergence result is plotted in Figure 13.

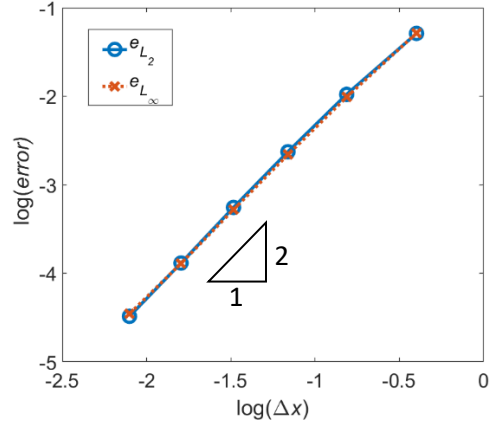

Figure 13. Convergence of the FC-RKPM solution to the exact solution for the 3D Poisson problem.

Again, the traditional RKPM solution of the same problem with the same nodes on Ω yields identical results.

### 6.2. Performance

To show the ultra-high-speed performance of the new method, we compare the run-times of the FC-RKPM with those from the traditional RKPM based on direct summations and matrix-vector operations. To have a fair comparison, both methods are coded in MATLAB and executed on the same computer. Each pair of simulations use identical nodes on Ω. Consequently, given the same problem and nodal spacing, the total number of nodes in FC-RKPM are slightly higher than the traditional RKPM due to the extension of Ω to $\mathbb{T}$.

For the traditional method, the stiffness matrix is constructed and stored, and the matrix-vector product is used to compute $\mathbf{f}^{int} = \mathbf{Kd}$ as needed. We used the most efficient serial implementation of the traditional method which is to have nested loops over neighbors only for computing non-zero elements of **K**, leading to a complexity of $O(NM^2)$ instead of $O(N^3)$ that one would get for a full-scale triple nested loop. For matrix storage and the matrix-vector products the sparse operations in MATLAB are used which resulted in the memory allocation of $O(NM)$ and the complexity of $O(NM)$ for the product. This is to maximize the efficiency and to minimize the memory allocation with the traditional method.

The tests were conducted on a Dell Precision 3650 desktop computer with Intel(R) Core(TM) i7-10700 CPU @ 2.90GHz and an Installed RAM of 32 GB. MATLAB R2021b was used for executing the codes. We used only one processor to compare the performance for serial implementation.

In the following, we compare the performance of the two methods for various terms in the analysis. We study the influence of the total node number ($N$), number of neighbors ($M$) and also the degree of monomial basis functions ($n$) on the CPU run-times and on the memory allocation for each method.

#### 6.2.1. Number of nodes

Here we study the influence of the node numbers/degrees of freedom on the performance of the two methods.

We begin the comparison by the most important term, that is, computing $\mathbf{f}^{int}$ in 3D where the new method offers the most efficiency gain. Table 1 show the results for the 3D problem described in Section 6.1.3. For the traditional method, we show the time needed to compute the stiffness matrix **K** and the **Kd** product in two separate columns in order to provide more informative data. Note that in these simulations the normalized support size is: $\widetilde{a_x} = \widetilde{a_y} = \widetilde{a_z} = 1.5$ which leads to $M = 27$ (number of neighbors) for interior nodes. Also, the linear monomial basis is used ($n = 1$) resulting in $s = 4$.

Table 1. CPU run-times for $\mathbf{f}^{\text{int}}$ computed by FC-RKPM and the traditional RKPM on a single processor.

| Number of Nodes on | FC-RKPM | RKPM | |
|---|---|---|---|
| | | K | Kd |
| 343 | $0.3\ ms$ | $6.54\ s$ | $78.5\ \mu s$ |
| 3,375 | $2.1\ ms$ | $96.6\ s$ | $0.30\ ms$ |
| 29,791 | $9.9\ ms$ | $17.5\ min$ | $2.87\ ms$ |
| 250,047 | $0.13\ s$ | $2.6\ hrs$ | $25.6\ ms$ |
| 2,048,383 | $1.22\ s$ | $22.4\ hrs$ | $0.26\ s$ |
| 16,581,375 | $10.3\ s$ | $\sim 7.5\ dys^*$ | $\sim 2.1\ s^*$ |

*estimated (insufficient memory)*

We observe a speed up of about 70,000 for computing $\mathbf{f}^{\text{int}}$ in 3D for $\widetilde{a_x} = \widetilde{a_y} = \widetilde{a_z} = 1.5$, and $n = 1$.

Note that the traditional solver computes **K** only once, and then uses the product with **d** to evaluate $\mathbf{f}^{\text{int}}$ every time afterwards. As a result, the computational time for the solver to finish would comprise of one **K** and $N_t$ (number of iterations in the implicit or time steps an explicit solver) times the **Kd** evaluation time. In FC-RKPM however, the run-time is the same for every evaluation of $\mathbf{f}^{\text{int}}$.

Table 2 shows the run time comparison for evaluating the other terms: $\mathbf{f}^r$, $\mathbf{u}^h$, and $\mathbf{M}$ & $\mathbf{M}^{-1}$.

Table 2. CPU run-times for $\mathbf{f}^r$, $\mathbf{u}^h$, and $\mathbf{M}$ & $\mathbf{M}^{-1}$ computed by FC-RKPM and the traditional RKPM.

| Number of Nodes on | $\mathbf{f}^r$ | | $\mathbf{u}^h$ | | $\mathbf{M}$ & $\mathbf{M}^{-1}$ | |
|---|---|---|---|---|---|---|
| | FC-RKPM | RKPM | FC-RKPM | RKPM | FC-RKPM | RKPM |
| 343 | $19.6\ ms$ | $23.5\ ms$ | $0.89\ ms$ | $18.7\ ms$ | $8.60\ ms$ | $33.4\ ms$ |
| 3,375 | $1.94\ ms$ | $0.17\ s$ | $2.03\ ms$ | $0.16\ s$ | $30.8\ ms$ | $0.30\ s$ |
| 29,791 | $4.45\ ms$ | $1.60\ s$ | $4.74\ ms$ | $1.45\ s$ | $0.23\ s$ | $2.41\ s$ |
| 250,047 | $42.2\ ms$ | $14.23\ s$ | $59.0\ ms$ | $12.6\ s$ | $1.69\ s$ | $21.3\ s$ |
| 2,048,383 | $0.33\ s$ | $1.9\ min$ | $0.46\ s$ | $1.78\ min$ | $14.1\ s$ | $2.93\ min$ |
| 16,581,375 | $2.86\ s$ | $\sim 15\ min^*$ | $3.72\ s$ | $\sim 14\ min^*$ | $2.1\ min$ | $\sim 24\ min^*$ |

*estimated (insufficient memory)*

We observe speed ups of about 200 to 300 for $\mathbf{f}^r$ and $\mathbf{u}^h$. For $\mathbf{M}$ & $\mathbf{M}^{-1}$, the speed up is about 12 at most, which is not as large as the speedups achieved for other terms. The reason is that computing $\mathbf{M}$ is the part that is significantly transformed by FC-RKPM; but the inversion operation is the same in both methods.

Note that in this 3D example, the number of neighbors is about the minimum possible (only nearest neighbors are covered by the shape functions). If the support size is selected to include more neighbors, the efficiency gains would become even higher. Because the complexity of FC-RKPM depends on $N$ only, but the performance of the conventional method would severely affected by larger $M$. This is investigated in the next section.

*Remark:* In explicit RKPM implantations where $\mathbf{f}^{\text{int}}$ is evaluated from $\mathbf{d}$ without storing the stiffness matrix (such as in [62]), the efficiency gains are expected to be about twice speed up as FC-RKPM offers in evaluating $\mathbf{u}^h$ (about 500). Because $\mathbf{f}^{\text{int}}$ requires two convolutions, not one.

*Remark:* Note that FC-RKPM uses more nodes in total (because of the gap between $\Omega$ and the box $\mathbb{T}$). As a result, geometries that are long and curved result in large gaps. Consequently, FC-RKPM would have nodes that are several times more compared to what is needed by the traditional RKPM for the same nodal spacing. Nevertheless, the efficiency gains by FC-RKPM are too high to be compromised by the extra nodes in the gap. Below, the memory allocation needed by the two methods to solve the 3D static problem is listed in Table 3.

Table 3. memory allocation in FC-RKPM and the traditional RKPM in 3D example with $\tilde{a} = 1.5$ and $n = 1$.

| Number of Nodes on | FC-RKPM | RKPM |
|---|---|---|
| 343 | 733 KB | 1.06 MB |
| 3,375 | 5.72 MB | 13.8 MB |
| 29,791 | 45.8 MB | 137 MB |
| 250,047 | 366 MB | 1.19 GB |
| 2,048,383 | 2.86 GB | 10.0 GB |
| 16,581,375 | 22.9 GB | ~82 GB* |

*estimated (insufficient memory)*

As observed for memory required in FC-RKPM for this example is about one-third the amount needed by the traditional method. As mentioned earlier, efficiency gains in memory allocation are much higher when more neighbors are covered in the shape function. This is demonstrated in the following section.

### 6.2.2. Number of neighbors and degree of monomial basis

In this part we study the influence of the number of neighbors, i.e., nodes that are covered by the shape functions ($M$), and also the influence of the monomial basis degree ($n$) on the CPU run-times and the memory allocation in the traditional and the new RKPM.

To this aim, we solved the 3D example described by Eq. (106), using a fixed number of nodes on $\Omega$, but with varying $M$ and $n$. The cubic $\Omega$ was discretized with $N_x = N_y = N_z = 20$ ($N = 8,000$ nodes on $\Omega$). Six simulations were performed by each method: three with $n = 1$ and $\tilde{a} = 1.5, 2.5, 3.5$; and three with $n = 2$ and $\tilde{a} = 2.5, 3.5, 4.5$. Since the number of nodes on $\Omega$ was fixed, we used Eq. (66) to find the extended domain for FC-RKPM. Table 4 and Table 5 show the CPU run-times for these simulations.

Table 4. CPU run-times for $\mathbf{f}^{\text{int}}$, $\mathbf{f}^{\text{r}}$, $\mathbf{u}^h$, and $\mathbf{M}$ & $\mathbf{M}^{-1}$ computed by FC-RKPM and the traditional RKPM using 8,000 nodes and $n = 1$.

| $\tilde{a}$ | $M$ | $\mathbf{f}^{\text{int}}$ | | $\mathbf{f}^{\text{r}}$ | | $\mathbf{u}^h$ | | $\mathbf{M}$ & $\mathbf{M}^{-1}$ | |
|---|---|---|---|---|---|---|---|---|---|
| | | FC-RKPM | RKPM | FC-RKPM | RKPM | FC-RKPM | RKPM | FC-RKPM | RKPM |
| 1.5 | 27 | 4.23 $ms$ | 4.58 $min$ | 2.68 $ms$ | 427 $ms$ | 2.45 $ms$ | 435 $ms$ | 72.1 $ms$ | 691 $ms$ |
| 2.5 | 125 | 6.66 $ms$ | 1.61 $hrs$ | 2.79 $ms$ | 1.83 $s$ | 2.50 $ms$ | 1.71 $s$ | 88.4 $ms$ | 2.54 $s$ |
| 3.5 | 343 | 7.71 $ms$ | 10.3 $hrs$ | 3.61 $ms$ | 4.51 $s$ | 3.93 $ms$ | 4.29 $s$ | 92.8 $ms$ | 6.17 $s$ |

Table 5. CPU run-times for $\mathbf{f}^{\text{int}}$, $\mathbf{f}^{\text{r}}$, $\mathbf{u}^h$, and $\mathbf{M}$ & $\mathbf{M}^{-1}$ computed by FC-RKPM and the traditional RKPM using 8,000 nodes and $n = 2$.

| $\tilde{a}$ | $M$ | $\mathbf{f}^{\text{int}}$ | | $\mathbf{f}^{\text{r}}$ | | $\mathbf{u}^h$ | | $\mathbf{M}$ & $\mathbf{M}^{-1}$ | |
|---|---|---|---|---|---|---|---|---|---|
| | | FC-RKPM | RKPM | FC-RKPM | RKPM | FC-RKPM | RKPM | FC-RKPM | RKPM |
| 2.5 | 125 | $10.1\ ms$ | $1.85\ hrs$ | $3.86\ ms$ | $1.95\ s$ | $4.65\ ms$ | $1.97\ s$ | $211\ ms$ | $2.93\ s$ |
| 3.5 | 343 | $29.4\ ms$ | $11.5\ hrs$ | $7.48\ ms$ | $4.93\ s$ | $9.25\ ms$ | $4.72\ s$ | $313\ ms$ | $7.11\ s$ |
| 4.5 | 729 | $16.6\ ms$ | $1.70\ dys$ | $5.97\ ms$ | $10.2\ s$ | $6.77\ ms$ | $9.41\ s$ | $259\ ms$ | $14.6\ s$ |

We observe that as $M$ increases for traditional method, the CPU time increases significantly as well; especially for $\mathbf{f}^{\text{int}}$. The FC-RKPM, however, is hardly affected by the increase in $M$. The reason is that $\mathbf{f}^{\text{int}}$ computational complexity in RKPM scales with $O(NM^2)$ and $\mathbf{f}^{\text{r}}$, $\mathbf{u}^h$ and $\mathbf{M}_{pq}$ scale with $O(NM)$, but in the case of FC-RKPM all terms scale with $O(N\log_2 N)$, meaning that *FC-RKPM is independent of the support size and numbers of neighbors*. The minor increase in FC-RKPM run-times is due to the increase in the extensions that depends on $\tilde{a}$ (see Eq. (66)), leading to slight increase in total number of nodes for FC-RKPM. Here, we used only 8000 nodes; in the case of large problems the influence of extension vanishes.

In Table 5, the run-times of FC-RKPM is surprisingly higher for $\tilde{a} = 3.5$ compared to 4.5. Here is the reason: since we wanted to study the influence of $\tilde{a}$ we kept the nodes on $\Omega$ fixed, i.e., we chose the nodal spacing first and then we computed the extensions to form the periodic box. This led to total node numbers that are not powers of two, and therefore FFT, operations are not carried out optimally and their run-time is affected depending on the constituent prime numbers of $N_x$, $N_y$ and $N_z$. However, this is not a concern as the speedups compared to the traditional RKPM are too high to be affected by it. Below are the tables with speedups computed from Table 4 and Table 5:

Table 6. Speedups for using FC-RKPM versus traditional RKPM in 3D and with $n = 1$.

| $\tilde{a}$ | $M$ | $\mathbf{f}^{\text{int}}$ | $\mathbf{f}^{\text{r}}$ | $\mathbf{u}^h$ | $\mathbf{M}$ & $\mathbf{M}^{-1}$ |
|---|---|---|---|---|---|
| 1.5 | 27 | 6.5e+2 | 1.6e+2 | 1.8e+2 | 9.6 |
| 2.5 | 125 | 8.7e+5 | 6.6e+2 | 6.8e+2 | 2.9e+1 |
| 3.5 | 343 | 4.8e+6 | 1.2e+3 | 1.1e+3 | 6.6e+1 |

Table 7. Speedups for using FC-RKPM versus traditional RKPM in 3D and with $n = 2$.

| $\tilde{a}$ | $M$ | $\mathbf{f}^{\text{int}}$ | $\mathbf{f}^{\text{r}}$ | $\mathbf{u}^h$ | $\mathbf{M}$ & $\mathbf{M}^{-1}$ |
|---|---|---|---|---|---|
| 2.5 | 125 | 6.6e+5 | 5.1e+2 | 4.2e+2 | 1.4e+1 |
| 3.5 | 343 | 1.4e+6 | 6.6e+2 | 5.1e+2 | 2.3e+1 |
| 4.5 | 729 | 8.8e+6 | 1.7e+3 | 1.4e+3 | 5.6e+1 |

The first row of Table 6 confirms the speedups we previously observed in Table 1 and Table 2. Comparing second and third rows of Table 6 with the first and second rows of Table 7, respectively, show that efficiency gains are slightly lower for larger $n$. The reason is the increase in $s$ from 4 (when $n = 1$) to 10 (when $n = 2$), increases the computational time in both methods, and therefore, the speedup being the ratio decreases.

We observe that as $M$ increases, the speedups for $\mathbf{M}$ & $\mathbf{M}^{-1}$ becomes higher. This is because the run-time fraction associated with computing $\mathbf{M}$ is significantly affected in the traditional method, whereas the FC-RKPM is not affected by the change in neighbor numbers.

If the number of neighbors increases, the required memory for the traditional RKPM will increase proportionally due to allocation of $N \times M$ arrays. In contrary, the memory needed in FC-RKPM will not be affected since the arrays have $N$ total elements at most. This is shown in Table 8 and Table 9.

Table 8. memory allocation in FC-RKPM and the traditional RKPM for the described 3D example with 8,000 nodes and $n = 1$.

| $\tilde{a}$ | $M$ | **FC-RKPM** | **RKPM** |
|---|---|---|---|
| 1.5 | 27 | 12.9 MB | 34.7 MB |
| 2.5 | 125 | 14.9 MB | 159 MB |
| 3.5 | 343 | 17.0 MB | 398 MB |

Table 9. memory allocation in FC-RKPM and the traditional RKPM for the described 3D example with 8,000 nodes and $n = 2$.

| $\tilde{a}$ | $M$ | **FC-RKPM** | **RKPM** |
|---|---|---|---|
| 2.5 | 125 | 43.6 MB | 170 MB |
| 3.5 | 343 | 49.9 MB | 409 MB |
| 4.5 | 729 | 56.6 MB | 748 MB |

The minor increase in the memory required by FC-RKPM is due to the slight increase of the extended domain, which becomes negligible for large problems.

*Remark*: in this work we only studied the efficiency gains in 3D. The speedups in 2D are expected to be lower than 3D but still significantly high.

### 6.2.3. High-Performance Computing

Another advantage of FC-RKPM is its low barrier in utilizing high-performance computing (HPC). While the fast convolution approach offers an ultra-high-speed serial implementation of RKPM, extra boosts in speed can be easily achieved with employing parallel FFT operations. Due to the wide applications of FFT, robust efficient parallel FFT/iFFT libraries already exist and can be called in FC-RKPM solvers with no additional effort. For example, FFT functions in MATLAB are capable of performing the operations in multithreaded fashion or on GPU if desired. For multithreading one just needs to specify the number of processors, and for GPU, arrays' type should be converted to "gpu array". Table 10 compares a serial and a GPU-based FC-RKPM simulation of the 3D problem using 2,048,383 nodes on Ω, $\tilde{a} = 3.5$, and $n = 1$. The GPU computations are carried out on a NVIDIA Quadro P2200 GPU with 5 GB memory, on the same desktop computer where the rest of the simulations were performed.

Table 10. Run-times for $\mathbf{f}^{\text{int}}$, $\mathbf{f}^{\text{r}}$, $\mathbf{u}^h$, and $\mathbf{M}$ & $\mathbf{M}^{-1}$ computed by FC-RKPM on a single CPU and on GPU using $N = 2,048,383$; $M = 343$, and $n = 1$.

|  | $\mathbf{f}^{\text{int}}$ | $\mathbf{f}^{\text{r}}$ | $\mathbf{u}^h$ | $\mathbf{M}$ & $\mathbf{M}^{-1}$ |
|---|---|---|---|---|
| 1 CPU | $1.18\ s$ | $0.29\ s$ | $0.48\ s$ | $14.1\ s$ |
| 1 GPU | $0.305\ s$ | $0.16\ s$ | $0.10\ s$ | $11.9\ s$ |

As observed the simulations can ran several times faster with minor modification of the serial code. As before, $\mathbf{M}$ & $\mathbf{M}^{-1}$ are least affected term since the inversion operation is serial in both rows and obtaining $\mathbf{M}$ is the part where FFT is used, and hence, benefits from GPU.

Based on the data on Table 1, Table 2, Table 4, and the complexity of the traditional method, this simulation takes about 4 month using the serial traditional RKPM. Indeed, parallel versions of the traditional method would improve performance, but requires high lever programming skills and immense computational power.

## 7. Conclusions

In this study, we introduced the fast-convolving reproducing kernel particle method (FC-RKPM). In this method, the governing equations are first discretized using reproducing kernel (RK) approximation. Then, the discrete system is expressed in terms of convolution sums, meaning that all summations arising from RK approximation or numerical quadrature should be expressed in convolutional forms. The convolutions are then efficiently computed using fast Fourier transform (FFT) and inverse FFT operations. Since all summations are computed in the Fourier space, no storage and looping over neighbors are required. The major cost is only the FFT operations which have the complexity of $O(N\log_2 N)$ with $N$ being the total number of nodes. In contrast, the most efficient traditional RKPM solvers require looping over neighbors for every node, resulting in double and triple nested loops with the complexity of $O(NM)$ and $O(NM^2)$, where $M$ is the number of neighbors and depends on the support size and nodal spacing. Removal of the nested loops by exploiting the convolutional structures of RK-based discretizations and the Convolution theorem leads to speedups ranging from hundreds to millions of times depending on the problem type, discretization method, and RK parameters.

For FFT to be applicable for fast convolutions on bounded arbitrary shaped domains (non-periodic conditions), we extend the domain of interest to a periodic box and modify the RK shape functions, such that the solution to the non-periodic problems are obtained, while allowing to use FFT for convolutions due to periodicity of the extended domain. As a model problem, the method is described in detail for the Galerkin weak form of Poisson problem with RK approximation. Implementation details and algorithms are provided, and the method's application to time-dependent and nonlinear problems is also discussed.

The FC-RKPM solutions are verified against analytical solutions in 1D, 2D and 3D, and the method is shown to be optimally convergent. In fact, the FC-RKPM solutions are found to be identical to the solution obtained by the traditional RKPM. We have provided a comprehensive comparison between the performance of FC-RKPM and the traditional method in 3D. We studied the influence of the number of nodes, number of neighbors (the support size), and the degree of the monomial basis functions on the performance of the new and the traditional method. Our results show speedups from hundreds to tens of thousands for the minimum neighbor numbers, and speedups as high as millions when more neighbors are considered. Days of simulations are now possible within seconds, and years of impractical meshfree simulations can be conducted in a few days on a single processor. FC-RKPM has the lowest barrier for high-performance computing, since parallel FFT libraries can easily be called instead of the serial ones with almost no programming effort, leading to some extra speedups.

**Acknowledgements**

The authors would like to thank Dr. Jiarui Wang of Brown University for the helpful discussions on the traditional RK methods and their development history.

## Appendix A.    FC-RKPM for boundary integrals

This appendix describes the FC-RKPM for boundary integrals that may arise in the weak forms that use RK approximation. In part A.1, the convolutional form is discussed. The numerical implementation and corresponding algorithm follow next in part A.2.

### A.1.    Convolutional form

As mentioned in the Section 2.3, in the standard RKPM, the computational cost of boundary integrals like $\mathbf{f}^q$ in Eq. (28), is far less than the cost of computing volume integrals like $\mathbf{K}$ and $\mathbf{f}^r$ (Eqs. (26) and (27)). The reason is that the computation is only required for nodes close to the boundaries and that the quadrature is performed on a lower dimensional manifold (on the boundaries). However, if desired, one can use the fast convolution operation for computing boundary integrals as well. This could be desired for implementation convenience and uniformity within FC-RKPM codes. Here, we show how the convolutional form for a boundary term can be achieved. Take $\mathbf{f}^q$ for example. Substituting Eq. (23) in Eq. (28), and using Eq. (41) yields:

$$f_I^q = \int_\Gamma \Psi_I(\mathbf{x}) q(\mathbf{x})\,\mathrm{d}\Gamma = \int_\Gamma \chi(\mathbf{x})\chi(\mathbf{x}_I)[\boldsymbol{b}^0]^\mathrm{T}(\mathbf{x}) \boldsymbol{H}^a(\mathbf{x}-\mathbf{x}_I) q(\mathbf{x})\,\mathrm{d}\Gamma \tag{A.1}$$

Using DNI for quadrature one gets:

$$f_I^q = \chi_I \sum_{S=1}^{M_I^\Gamma} \chi_S [\boldsymbol{b}^0]_S^\mathrm{T} \boldsymbol{H}_{S-I}^a q_S A_S \tag{A.2}$$

Where $A_S$ is the quadrature weight, i.e., the area associated with the boundary node $\mathbf{x}_S$. $M_I^\Gamma$ is the number of nodes that are neighbors of $\mathbf{x}_I$ and locate on the boundary $\Gamma$.

Knowing that $q_S$ and $A_S = 0$ for all $\mathbf{x}_S \notin \Gamma$, one can extend the summation to include all the nodes in the domain $\mathbb{T}$:

$$f_I^q = \chi_I \sum_{S=1}^{N} \chi_S [\boldsymbol{b}^0]_S^T \boldsymbol{H}_{S-I}^a q_S A_S \quad (A.3)$$

Using indicial notation for the vector product, and using Eq. (55):

$$f_I^q = \chi_I \sum_{S=1}^{N} \chi_S q_S A_S (b_p^0)_S (\overline{H_p^a})_{I-S} = \chi_I (\chi q A b_p^0 \circledast \overline{H_p^a})_I \quad , \text{and } p = 1,2,\ldots,s \quad (A.4)$$

Note that we used $\mathbf{f}^q$ as an example here; the procedure is general and can be carried out for any other boundary integral.

*Remark:* Note that the computational cost of the new approach for boundary integrals is in the same order as the body terms $\mathbf{f}^{int}$ and $\mathbf{f}^r$ with FC-RKPM. Although the quadrature is accelerated via FFT, its summation domain is extended; therefore, the computational gain compared with the standard method is unknown, and should be investigated in future.

### A.2. Numerical implementation and algorithm

Having the necessary arrays computed and stored, the computer implementation for fast evaluation of $\mathbf{f}^q$ in Eq. (A.4) is presented as follows:

$$\mathbf{f}^q = \boldsymbol{\chi} \circ \mathbf{F}^{-1}[\mathbf{F}(\boldsymbol{\chi} \circ \mathbf{A} \circ \mathbf{b}_p^0 \circ \mathbf{q}) \circ \mathbf{F}(\overline{\mathbf{H}_p^a})] \quad (A.5)$$
, and $p = 1,2,\ldots,s$.

where

$$\mathbf{q} = [q_{ij}]_{N_1 \times N_2} = [q(\mathbf{x}_{ij})]_{N_1 \times N_2} \quad (A.6)$$
$$\mathbf{A} = [A_{ij}]_{N_1 \times N_2} = [A(\mathbf{x}_{ij})]_{N_1 \times N_2} \quad (A.7)$$

The subroutine for the boundary integral given by Eq. (A.5), is then:

---
Algorithm A.1.  Subroutine for the boundary integral given by Eq. (A.5): $\mathbf{f}^q$

---

Precomputed in initialization: $\mathbf{b}_p^0$; $\mathbf{G}_p^0 = \boldsymbol{\chi} \circ \mathbf{A} \circ \mathbf{b}_p^0$; $\widehat{\overline{\mathbf{H}_p^a}} = \mathbf{F}(\overline{\mathbf{H}_p^a})$; (for all $p = 1:s$)

Given $\mathbf{q}$ as input:

Set $\widehat{\mathbf{B}} = 0$
For $p = 1:s$
$\quad \widehat{\mathbf{B}} = \widehat{\mathbf{B}} + \mathbf{F}(\mathbf{G}_p^0 \circ \mathbf{q}) \circ \widehat{\overline{\mathbf{H}_p^a}}$
End

$\mathbf{f}^q = \boldsymbol{\chi} \circ \mathbf{F}^{-1}(\widehat{\mathbf{B}})$

---

### Appendix B.  FC-RKPM for nonlinear problems

This appendix describes the FC-RKPM for nonlinear problems. First the convolutional forms of the nonlinear terms are discussed. The numerical implementation and corresponding algorithms follow next.

### B.1. Convolutional forms

For a PDE with a nonlinear term $\mathcal{N}(u)$, the corresponding integral in the Galerkin weak form can be of form:

$$\int_\Omega w^h \mathcal{N}(u^h) \mathrm{d}\Omega \tag{B.1}$$

Using the vector forms in Eq. (18), one gets the following vector which we refer to as the *nonlinear force*:

$$\mathbf{f}^\mathcal{N} = \int_\Omega \mathbf{N}^\mathrm{T} \mathcal{N}(u^h) \mathrm{d}\Omega \tag{B.2}$$

or

$$f_I^\mathcal{N} = \int_\Omega \Psi_I(\mathbf{x}) \mathcal{N}_u(\mathbf{x}) \mathrm{d}\Omega = \int_\mathbb{T} \chi(\mathbf{x}) \chi(\mathbf{x}_I)[\mathbf{b}^0]^\mathrm{T}(\mathbf{x}) H^a(\mathbf{x} - \mathbf{x}_I) \mathcal{N}_u(\mathbf{x}) \mathrm{d}\mathbb{T} \tag{B.3}$$

where $\mathcal{N}_u(\mathbf{x}) = \mathcal{N}\left(u^h(\mathbf{x})\right)$. Using DNI for quadrature and indicial notation of the vector product yields:

$$f_I^\mathcal{N} = \sum_{S=1}^N \chi(\mathbf{x}_S) \chi(\mathbf{x}_I) b_p^0(\mathbf{x}_S) H_p^a(\mathbf{x}_S - \mathbf{x}_I) \mathcal{N}_u(\mathbf{x}_S) V(\mathbf{x}_S) \tag{B.4}$$

$$= \chi_I \sum_{S=1}^N \chi_S (b_p^0)_S \mathcal{N}_{u_S} V_S (\overline{H_p^a})_{I-S} = \chi_I \{[\chi b_p^0 \mathcal{N}_u V] \circledast \overline{H_p^a}\}_I$$

and $p = 1,2,\ldots,s$

$\mathcal{N}_u(\mathbf{x}_I)$ is computed given the approximated function $u^h$ at the latest time step (in explicit solvers), or the latest iteration (in iterative implicit solvers).

Another term that is likely to arise in the weak form of nonlinear problems is:

$$\int_\Omega \nabla w^h \cdot \mathcal{N}(u^h) \mathrm{d}\Omega \tag{B.5}$$

where $\mathcal{N}(u^h)$ has the same dimension as $\nabla w^h$ and is a nonlinear function of $u^h$. Using Eq. (20), the resulting vector from becomes:

$$\mathbf{f}^\mathcal{N} = \int_\Omega \mathbf{B}^\mathrm{T} \mathcal{N}(u^h) \mathrm{d}\Omega \tag{B.6}$$

$\mathbf{B}$ and $\mathcal{N}$ then depend on the field variable being scalar or vector, and also, on the spatial dimension of the problem. For a scalar problem in 2D, for example, one gets:

$$f_I^\mathcal{N} = \int_\Omega \left( \Psi_I^{\nabla_x}(\mathbf{x}) \mathcal{N}_u^x(\mathbf{x}) + \Psi_I^{\nabla_y}(\mathbf{x}) \mathcal{N}_u^y(\mathbf{x}) \right) \mathrm{d}\Omega \tag{B.7}$$

$$= \int_\mathbb{T} \chi(\mathbf{x}) \chi(\mathbf{x}_I) \{[\mathbf{b}^x]^\mathrm{T}(\mathbf{x}) \mathcal{N}_u^x(\mathbf{x}) + [\mathbf{b}^y]^\mathrm{T}(\mathbf{x}) \mathcal{N}_u^y(\mathbf{x})\} H^a(\mathbf{x} - \mathbf{x}_I) \mathrm{d}\mathbb{T}$$

Using DNI for quadrature and indicial notation:

$$f_I^{\mathcal{N}} = \sum_{S=1}^{N} \chi(x_S)\chi(x_I)\big[b_p^x(x_S)\mathcal{N}_u^x(x_S) + b_p^y(x_S)\mathcal{N}_u^y(x_S)\big]H_p^a(x-x_I)V(x_S) \quad \text{(B.8)}$$

$$= \chi_I \sum_{S=1}^{N} \chi_S \left[(b_p^x)_S \mathcal{N}_u^x{}_S + (b_p^y)_S \mathcal{N}_u^y{}_S\right] V_S \big(\overline{H_p^a}\big)_{I-S}$$

$$= \chi_I \big\{\big[\chi(b_p^x \mathcal{N}_u^x + b_p^y \mathcal{N}_u^y)V\big] \circledast \overline{H_p^a}\big\}_I;$$

and $p = 1,2,\ldots,s$

Another approach to tackle nonlinear problems is, of course, to linearize the system [63] which would then lead to a similar structure as the one in the linear internal force ($\mathbf{f}^{int} = \mathbf{Kd}$) which is described in Section 4.2.2, and can be tailored accordingly.

### B.2. Numerical implementation and algorithms

Having the necessary arrays computed and stored, the computer implementation for fast evaluation of the $\mathbf{f}^{\mathcal{N}}$ in Eq. (B.4) becomes:

$$\mathbf{f}^{\mathcal{N}} = \boldsymbol{\chi} \circ \mathbf{F}^{-1}\big[\mathbf{F}(\boldsymbol{\chi} \circ \mathbf{V} \circ \mathbf{b}_p^0 \circ \mathbf{N}_u) \circ \mathbf{F}(\overline{H_p^a})\big], \text{ and } p = 1,2,\ldots,s. \quad \text{(B.9)}$$

The $\mathbf{f}^{\mathcal{N}}$ given by Eq. (B.8) becomes:

$$\mathbf{f}^{\mathcal{N}} = \boldsymbol{\chi} \circ \mathbf{F}^{-1}\big\{\mathbf{F}\big[\boldsymbol{\chi} \circ \mathbf{V} \circ (\mathbf{b}_p^x \circ \mathbf{N}_u^x + \mathbf{b}_p^y \circ \mathbf{N}_u^y)\big] \circ \mathbf{F}(\overline{H_p^a})\big\}, \text{ and } p = 1,2,\ldots,s \quad \text{(B.10)}$$

Having $\mathbf{u}^h$ evaluated from Algorithm 3, a nonlinear internal force in the form of Eq. (B.9) is computed from the subroutine in Algorithm B.1.

---

Algorithm B.1.  Subroutine for evaluation of a nonlinear internal force density in the form of Eq. (B.9): $\mathbf{f}^{\mathcal{N}}(\mathbf{d})$

---

Precomputed in initialization: $\mathbf{b}_p^x; \mathbf{b}_p^y; \mathbf{C}_p^x = \boldsymbol{\chi} \circ \mathbf{V} \circ \mathbf{b}_p^x; \mathbf{C}_p^y = \boldsymbol{\chi} \circ \mathbf{V} \circ \mathbf{b}_p^y; \widehat{\overline{H_p^a}} = \mathbf{F}(\overline{H_p^a});$ (for all $p = 1{:}s$),

Given $\mathbf{d}$

Evaluate $\mathbf{u}^h$ from Algorithm 3
Compute $\mathbf{N}_u^x$ and $\mathbf{N}_u^y$ using the given nonlinear functions of $\mathbf{u}^h$

Set $\widehat{\mathbf{B}} = 0$
For $p = 1{:}s$
 $\widehat{\mathbf{B}} = \widehat{\mathbf{B}} + \mathbf{F}(\mathbf{C}_p^x \circ \mathbf{N}_u^x + \mathbf{C}_p^y \circ \mathbf{N}_u^y) \circ \widehat{\overline{H_p^a}}$
End

$\mathbf{f}^{\mathcal{N}} = \boldsymbol{\chi} \circ \mathbf{F}^{-1}(\widehat{\mathbf{B}})$

---

If the nonlinear internal force is of the from in Eq. (B.10), then the algorithms become:

| Algorithm B.2. | Subroutine for evaluation of a nonlinear internal force density in the form of Eq. (B.10): $\mathbf{f}^{\mathcal{N}}(\mathbf{d})$ |
|---|---|

Precomputed in initialization: $\mathbf{b}_p^0$; $\mathbf{C}_p^0 = \boldsymbol{\chi} \circ \mathbf{V} \circ \mathbf{b}_p^0$; $\widehat{\overline{\mathbf{H}_p^a}} = \mathbf{F}(\overline{\mathbf{H}_p^a})$; (for all $p = 1:s$),

Given $\mathbf{d}$;

Evaluate $\mathbf{u}^h$ from Algorithm 3.
Compute $\mathbf{N}_u$ using the given nonlinear function of $\mathbf{u}^h$
Set $\widehat{\mathbf{B}} = 0$
For $p = 1:s$
$\quad \widehat{\mathbf{B}} = \widehat{\mathbf{B}} + \mathbf{F}(\mathbf{C}_p^0 \circ \mathbf{N}_u) \circ \widehat{\overline{\mathbf{H}_p^a}}$
End

$\mathbf{f}^{\mathcal{N}} = \boldsymbol{\chi} \circ \mathbf{F}^{-1}(\widehat{\mathbf{B}})$

The nonlinear terms require $2(s + 1)$ FFT/iFFT operations (including the evaluation of $\mathbf{u}^h$).

### B.3. Nonlinear solvers

For nonlinear problems, indeed various approaches can be adopted. One general approach is to use a nonlinear conjugate gradient method to solve the system for $\mathbf{d}$, i.e., to minimize the residual $\mathbf{R}(\mathbf{d})$ which is the discretized weak form of the PDE. Having a nonlinear term $\mathbf{f}^{\mathcal{N}}(\mathbf{d})$ like the ones given by Eq. (B.9) or Eq. (B.10) in $\mathbf{R}(\mathbf{d}) = 0$, nonlinear CG algorithms use successive evaluations of $\mathbf{f}^{\mathcal{N}}(\tilde{\mathbf{d}})$, with an updating $\tilde{\mathbf{d}}$ to find the minimizer $\mathbf{d}$. Therefore, the subroutines in Algorithm B.1 and Algorithm B.2 with the appropriate input $\tilde{\mathbf{d}}$ can be called in nonlinear CG iterations. FC-RKPM solver for nonlinear static is provided in Algorithm B.3.

| Algorithm B.3. | FC-RKPM solver for nonlinear static (time-independent) analysis. |
|---|---|

***Read inputs:***
*Discrete domain & geometry inputs:* $\mathbf{X}, \mathbf{Y}, \mathbf{V}, \boldsymbol{\chi}, \boldsymbol{\chi}_{\Gamma_g}$;
*Physical inputs:* material properties, boundary conditions' info, body force/heat source $\mathbf{r}(\mathbf{X}, \mathbf{Y})$;
*RKPM inputs:* $s, n, a_x, a_y$;
*Solver inputs:* e.g., tolerance
***Initialize:***
Define functions: $\Phi(\mathbf{X}, \mathbf{Y})$, $\mathbb{H}_p(\mathbf{X}, \mathbf{Y})$ for $p = 1, \dots, s$
Compute and store: $\mathbf{H}_p, \mathbf{H}_p^a, \overline{\mathbf{H}_p^a}$ (for $p = 1, \dots, s$) from Eqs. (77)(78) to (79); $\mathbf{b}_p^0, \mathbf{b}_p^x, \mathbf{b}_p^y$ from Algorithm 1; and other arrays needed in the subroutines to be called
Compute external force: $\mathbf{f}^r$ (Algorithm 4)
Compute boundary terms: e.g., $\mathbf{f}^q$ from Algorithm A.1
Set $\mathbf{d} = 0$ initially
***Solver:*** nonlinear conjugate gradient
While error > tolerance: compute $\mathbf{f}^{int}$ (Algorithm B.1 or Algorithm B.2) as needed and update $\mathbf{d}$
***Output:***
Compute $\mathbf{u}^h$ (Algorithm 3)
Write outputs

## Appendix C. FC-RKPM for time-dependent problems

In RKPM discretization of time dependent problems such as transient diffusion or equation of motion, the Galerkin weak form results in a term associated with the time derivative of the field variable. Take diffusion equation, for example, with the field variable $u(x,t)$ being a function of space and time:

$$\begin{cases} \dot{u} = \nu \nabla^2 u + r, & \text{on } \Omega \ (t > 0) \\ u = u_0, & \text{on } \Omega \ (t = 0) \\ u = g, & \text{on } \Gamma_g (t \geq 0) \\ \nabla u \cdot \boldsymbol{n} = q, & \text{on } \Gamma_q (t \geq 0) \end{cases} \tag{C.1}$$

where $\dot{u}$ denotes $\frac{\partial u(x,t)}{\partial t}$, $\nu$ is the diffusion coefficient, and $u_0$ is the initial condition. Let $\nu = 1$ for simplicity. The weak form is then:

$$\int_\Omega w^h \dot{u}^h \, d\Omega = -\int_\Omega \nabla w^h \cdot \nabla u^h \, d\Omega + \int_\Omega w^h r \, d\Omega + \int_{\Gamma_q} w^h q \, d\Gamma \tag{C.2}$$

Where:

$$\dot{u}^h(\boldsymbol{x},t) = \sum_{I=1}^N \Psi_I(\boldsymbol{x}) \dot{d}_I \tag{C.3}$$

Using Eqs.(17) to (24), Eq. (C.2) can be express as:

$$\mathbf{M}\dot{\mathbf{d}} + \mathbf{K}\mathbf{d} = \mathbf{f}^r + \mathbf{f}^q \tag{C.4}$$

where:

$$\mathbf{M} = \int_\Omega \mathbf{N}^T \mathbf{N} \, d\Omega \tag{C.5}$$

$$\dot{\mathbf{d}} = [\dot{d}_1, \dot{d}_2, \dots, \dot{d}_N]^T \tag{C.6}$$

### C.1. Convolutional forms

Similar to the wholistic approach used for reaching the convolutional form for $\mathbf{f}^{int} = \mathbf{K}\mathbf{d}$, one can find the convolutional structure for $\mathbf{f}^m = \mathbf{M}\dot{\mathbf{d}}$:

$$f_I^m = \sum_{J=1}^N M_{IJ} \dot{d}_J = \sum_{J=1}^N \left( \int_\mathbb{T} \Psi_I(\boldsymbol{x}) \Psi_J(\boldsymbol{x}) \, d\mathbb{T} \right) \dot{d}_J = \int_\mathbb{T} \left[ \sum_{J=1}^N \Psi_I(\boldsymbol{x}) \Psi_J(\boldsymbol{x}) \dot{d}_J \right] d\mathbb{T}$$

$$= \int_\mathbb{T} \left[ \sum_J \chi(\boldsymbol{x}) \chi(\boldsymbol{x}_I) [\boldsymbol{b}^0]^T(\boldsymbol{x}) H^a(\boldsymbol{x} - \boldsymbol{x}_I) \chi(\boldsymbol{x}) \chi(\boldsymbol{x}_J) [\boldsymbol{b}^0]^T(\boldsymbol{x}) H^a(\boldsymbol{x} - \boldsymbol{x}_J) \dot{d}_J \right] d\mathbb{T} \tag{C.7}$$

Using indicial notation for vector products and DNI for quadrature, and following the procedure shown in Eqs. (51) to (57) yields:

$$f_I^m = \chi_I \{ [\chi b_q^0 (\chi \dot{d} \circledast H_q^a) b_p^0 V] \circledast \overline{H_p^a} \}_I, \quad \text{and } p, q = 1, 2, \dots, s. \tag{C.8}$$

*Lumped mass:*

*Mass lumping* is a technique that, with little compromise in the accuracy, leads to more efficient time marching algorithms [64]. In this technique, the mass matrix **M**, is replaced with the *lumped mass* $\mathbf{M}^l$ which is a diagonal matrix, and the elements on the diagonal are obtained by summing over rows of **M**:

$$M_{II}^l = \sum_{P=1}^{N} M_{IP} = \sum_{P=1}^{N} \int_{\mathbb{T}} \Psi_I(x)\Psi_P(x)\, d\mathbb{T} = \int_{\mathbb{T}} \left[\sum_{P=1}^{N} \Psi_I(x)\Psi_P(x)\right] d\mathbb{T} \quad (C.9)$$

$$= \int_{\mathbb{T}} \Psi_I(x)\left[\sum_{P=1}^{N} \Psi_P(x)\right] d\mathbb{T} = \int_{\mathbb{T}} \Psi_I(x)\, d\mathbb{T}$$

Using Eqs. (41) and (55), and DNI for quadrature one gets:

$$M_{II}^l = \sum_{S=1}^{N} \chi_S \chi_I (b_p^0)_S (H_p^a)_{S-I} V_S = \chi_I \sum_{S=1}^{N} \chi_S V_S (b_p^0)_S \overline{(H_p^a)}_{I-S} \quad (C.10)$$

$$= \chi_I \left(\chi V b_p^0 \circledast \overline{H_p^a}\right)_I \quad , \text{and } p = 1,2,\ldots,s$$

*Remark:* The term associated with the time derivative in equation of motion is very similar. The only difference is that the order of time derivative for $d$ is two, i.e., $\mathbf{f}^m = \mathbf{M}\ddot{\mathbf{d}}$.

### C.2. Numerical implementation and algorithms

Having the necessary arrays computed and stored, the computer implementation for fast evaluation of $\mathbf{f}^m$ in Eq. (C.8) becomes:

$$\mathbf{f}^m = \boldsymbol{\chi} \circ \mathbf{F}^{-1}\left[\mathbf{F}(\boldsymbol{\chi} \circ \mathbf{V} \circ b_p^0 \circ \mathbf{A}^0) \circ \mathbf{F}\left(\overline{H_p^a}\right)\right], \text{and } p = 1,2,\ldots,s. \quad (C.11)$$

with

$$\mathbf{A}^0 = b_q^0 \circ \mathbf{F}^{-1}\left[\mathbf{F}(\boldsymbol{\chi} \circ \dot{\mathbf{d}}) \circ \mathbf{F}\left(H_q^a\right)\right], and\ q = 1,2,\ldots,s. \quad (C.12)$$

The lumped mass in Eq. (C.10) becomes:

$$\mathbf{M}^l = \boldsymbol{\chi} \circ \mathbf{F}^{-1}\left[\mathbf{F}(\boldsymbol{\chi} \circ \mathbf{V} \circ b_p^0) \circ \mathbf{F}\left(\overline{H_p^a}\right)\right], \text{and } p = 1,2,\ldots,s. \quad (C.13)$$

For the integral associated with the mass term given by Eq. (C.11), the following subroutine can be used.

| Algorithm C.1. | Subroutine for the mass term given by Eq.(C.11): $\mathbf{f}^m(\dot{\mathbf{d}})$ |
|---|---|

Precomputed in initialization: $\mathbf{b}_p^0; \mathbf{C}_p^0 = \boldsymbol{\chi} \circ \mathbf{V} \circ \mathbf{b}_p^0; \widehat{\mathbf{H}_p^a} = \mathbf{F}(\mathbf{H}_p^a); \widehat{\overline{\mathbf{H}_p^a}} = \mathbf{F}(\overline{\mathbf{H}_p^a});$ (for all $p = 1:s$)

Given $\dot{\mathbf{d}}$ as input:

$\widehat{\dot{\mathbf{d}}_\chi} = \mathbf{F}(\boldsymbol{\chi} \circ \dot{\mathbf{d}})$
Set $\mathbf{A}^0 = 0$;
For $p = 1:s$
$\quad \mathbf{D}_p = \mathbf{F}^{-1}\left[\widehat{\dot{\mathbf{d}}_\chi} \circ \widehat{\mathbf{H}_p^a}\right]$
$\quad \mathbf{A}^0 = \mathbf{A}^0 + \mathbf{b}_p^0 \circ \mathbf{D}_p$
End

$\widehat{\mathbf{B}} = 0$
For $p = 1:s$
$\quad \widehat{\mathbf{B}} = \widehat{\mathbf{B}} + \mathbf{F}(\mathbf{C}_p^0 \circ \mathbf{A}^0) \circ \widehat{\overline{\mathbf{H}_p^a}}$
End

$\mathbf{f}^m = \boldsymbol{\chi} \circ \mathbf{F}^{-1}(\widehat{\mathbf{B}})$

lumped mass can be evaluated using the following subroutine:

| Algorithm C.2. | Subroutine for the lumped mass given by Eq.(C.11) (C.13): $\mathbf{M}^l$ |
|---|---|

Precomputed in initialization: $\mathbf{b}_p^0; \mathbf{C}_p^0 = \boldsymbol{\chi} \circ \mathbf{V} \circ \mathbf{b}_p^0; \widehat{\overline{\mathbf{H}_p^a}} = \mathbf{F}(\overline{\mathbf{H}_p^a});$ (for all $p = 1:s$)

Set $\widehat{\mathbf{B}} = 0$
For $p = 1:s$
$\quad \widehat{\mathbf{B}} = \widehat{\mathbf{B}} + \mathbf{F}(\mathbf{C}_p^0) \circ \widehat{\overline{\mathbf{H}_p^a}}$
End

$\mathbf{M}^l = \boldsymbol{\chi} \circ \mathbf{F}^{-1}(\widehat{\mathbf{B}})$

### C.3. Explicit and implicit solvers

For solving time-dependent problems, time span is discretized into $N_t$ time steps, and a time integration scheme is used to update $\mathbf{d}$ at each time step.

If mass lumping is used, the following algorithm is suggested for explicit FC-RKPM analysis.

| Algorithm C.3. | FC-RKPM explicit analysis for time-dependent problems |

*Read inputs:*
*Discrete domain & geometry inputs:* $\mathbf{X}, \mathbf{Y}, \mathbf{V}, \boldsymbol{\chi}, \boldsymbol{\chi}_{\Gamma_g}$;
*Physical inputs:* material properties, Initial and boundary conditions' info, body force/heat source $\mathbf{r}(\mathbf{X}, \mathbf{Y}, t)$;
*RKPM inputs:* $s, n, a_x, a_y$;
*Solver inputs:* $N_t, t_{max}$
*Initialize:*
Define functions: $\Phi(\mathbf{X}, \mathbf{Y}), \mathbb{H}_p(\mathbf{X}, \mathbf{Y})$ for $p = 1, \ldots, s$
Compute and store: $\mathbf{H}_p, \mathbf{H}_p^a, \overline{\mathbf{H}_p^a}$ (for $p = 1, \ldots, s$) from Eqs. (77)(78) to (79); $\mathbf{b}_p^0, \mathbf{b}_p^x, \mathbf{b}_p^y$ from Algorithm 1; and other arrays needed in the subroutines to be called
Compute lumped mass: $\mathbf{M}^l$ (Algorithm C.2)
Compute external force*: $\mathbf{f}^r$ (Algorithm 4)
 (*If $\mathbf{f}^r$ is time-dependent, it should be inside the time loop)
Compute boundary terms*: e.g., $\mathbf{f}^q$ from Algorithm A.1
 (*If $\mathbf{f}^q$ is time-dependent, it should be inside the time loop)
Set $\mathbf{d} = \mathbf{d}(\mathbf{X}, \mathbf{Y}, t = 0)$ I.C.
*Solver:*
While $t < t_{max}$:
  $t = t + \Delta t$
  Compute $\mathbf{f}^{int}(\mathbf{d})$ (either Algorithm 2, Algorithm B.1, or Algorithm B.2)
  Compute $\dot{\mathbf{d}}$ (e.g., from Eq. (C.4) and using $\mathbf{M}^l$)
  Update $\mathbf{d}$ (using conventional explicit time integration methods)
*Output:*
Compute $\mathbf{u}^h$ (Algorithm 3)
Write outputs

If mass lumping is not used then, one needs to use a Krylov subspace iterative method to solve the system, e.g., Eq. (C.4). However, the cost would be similar to an implicit time integration which is more stable and allows for larger time steps. In this case, an implicit solver is preferred for time integration. The algorithm for implicit dynamic FC-RKPM analysis is given by Algorithm C.4.

| Algorithm C.4. | Solver structure for FC-RKPM implicit dynamic analysis. |

*Read inputs:*
⋮
*Solver inputs:* $N_t, t_{max}$, tolerance
*Initialize:*
⋮
*Solver:* (conventional implicit time integration methods)
While $t < t_{max}$:
  $t = t + \Delta t$
  While error > tolerance (conjugate gradient method)
    Compute $\mathbf{f}^{int}$ (either Algorithm 2, Algorithm B.1, or Algorithm B.2);
    and $\mathbf{f}^m$ (Algorithm C.1) as needed to find $\mathbf{d}$ at $t$
*Output:*
⋮

Note that the algorithms presented in this section are examples to demonstrate the general outline and structure of FC-RKPM solvers. Indeed, changes in the problems and discretization of different types of PDE's require modifications to the suggested algorithms.